\documentclass[preprint,10pt]{elsarticle}
\usepackage{graphicx}
\usepackage{amsmath}
\usepackage{subfigure}
\usepackage{longtable}
\usepackage{listings} 

\usepackage{amssymb}
\usepackage{amsthm}

\RequirePackage[top=1.5cm,bottom=1.5cm,left=2.cm,right=2.cm,includehead,includefoot]{geometry}
\RequirePackage[rgb,x11names]{xcolor}

\usepackage{indentfirst}

\journal{Journal of Computational Physics}
\begin{document}
\begin{frontmatter}



\title{A numerical study of two-phase flow with dynamic capillary pressure using an adaptive moving mesh method}

\author{Hong Zhang$^*$, \quad Paul Andries Zegeling}
\address{Department of Mathematics, Faculty of Science, Utrecht University, Budapestlaan 6, 3584CD Utrecht, The Netherlands }
\cortext[]{Corresponding author. \newline
E-mail addresses: H.Zhang4@uu.nl (H. Zhang), P.A.Zegeling@uu.nl (P. A. Zegeling)}
\begin{abstract}
Motivated by observations of saturation overshoot, this paper investigates numerical modeling of two-phase flow incorporating dynamic capillary pressure. The effects of the dynamic capillary coefficient, the infiltrating flux rate and the initial and boundary values are systematically studied using a travelling wave ansatz and efficient numerical methods. The travelling wave solutions may exhibit monotonic, non-monotonic or plateau-shaped behaviour. Special attention is paid to the non-monotonic profiles. The travelling wave results are confirmed by numerically solving the partial differential equation using an accurate adaptive moving mesh solver. Comparisons between the computed solutions using the Brooks-Corey model and the laboratory measurements of saturation overshoot verify the effectiveness of our approach.
\end{abstract}
\end{frontmatter}

\section{Introduction}
Since the proposition of the dynamic capillary concept \cite{barenblatt1971filtration, stauffer1978time, hassanizadeh1993thermodynamic}, the modified Buckley-Leverett (MBL) equation
\begin{align}
	\label{eqn:fullmbl0}
	\frac{\partial u}{\partial t} + \frac{\partial}{\partial x} F(u) = - \frac{\partial}{\partial x}[ H(u) \frac{\partial}{\partial x} (p_c(u) - \tau \frac{\partial u}{\partial t})],
\end{align}
which models the one-dimensional two-phase flow in porous media, has attracted considerable interest in hydrology and mathematics \cite{hassanizadeh2002dynamic, nieber2005dynamic, dicarlo2005modeling, van2007new, mikelic2010global, spayd2011buckley}. In the MBL equation, the third order mixed derivatives term represents the dynamic capillary pressure effect in the phase pressure difference. With the help of this term, Refs. \cite{dicarlo2005modeling, sander2008dynamic, chapwanya2010numerical, zhang2016mimetic} successfully captured the non-monotonic saturation profiles found by \cite{glass1989mechanism, selker1992fingered, liu1994formation, dicarlo2004experimental}. 

Different versions of the MBL equation have been studied from various points of view. Results on travelling wave solutions, global existence, phase plane analysis and uniqueness of weak solutions are given in \cite{cuesta2000infiltration, van2007new, van2013travelling, mikelic2010global, spayd2011buckley, cao2015uniqueness}. Ref. \cite{van2007new} shows that the travelling wave solutions of the MBL equation can be described by rarefaction wave, admissible Lax shock and undercompressive shock. In order to capture all these structures accurately, several numerical methods have been proposed in literature. A cell-centered finite difference method and a locally conservative Eulerian-Lagrangian method were proposed in Ref. \cite{peszynska2008numerical}, but it's mentioned that such methods may suffer from instabilities in convection-dominated cases and for large dynamic effects. Van Duijn et al. \cite{van2007new} applied a finite difference method which adopted a minmod slope limiter based on the first order upwind and Richtmyer's schemes. The solutions obtained by this method have good agreement with the travelling wave results. Wang and Kao \cite{wang2013central} extended the second and third order central schemes to capture the nonclassical solutions of the MBL equation. Kao et al. \cite{kao2015fast} split the MBL into a high-order linear equation and a nonlinear convective equation, and then integrated the linear equation with a pseudo-spectral method and the nonlinear equation with a Godunov-type central-upwind scheme. The computed solutions demonstrate that the higher-order spatial reconstruction using fifth-order WENO5 scheme gives more accurate numerical solutions. Zegeling \cite{zegeling2015adaptive} investigated the non-monotonic behaviour of a simple MBL equation with an adaptive moving mesh method, the result shows that for obtaining the same accuracy, the adaptive method needs around a factor of 4 fewer grid points than the uniform grid case.

As the moving mesh method has demonstrated outstanding advantages in tracking shocks or steep wave fronts of other two-phase flow equations \cite{doster2010numerical, hu2011simulating, dong2014adaptive}, in the present paper, we will study the solutions of the MBL equation using this method. To our best knowledge, the adaptive moving mesh method has not been applied to solve the full MBL equation which includes gravity and non-linear diffusion.


The rest of the paper is organized as follows. In Section 2, we outline the derivation of the two-phase flow equation and present some traveling wave analysis. Section 3 introduces a moving mesh method in terms of coordinate transformation. Numerical experiments are presented in Section 4 to show the effectiveness of the proposed method. Section 5 ends with the conclusion.

\section{Background}
In this section, we first derive the two-phase flow equation and then present some travelling wave results.

\subsection{The one-dimensional two-phase flow equation with dynamic capillary pressure term}\label{sec:equation}
Here we use the fractional flow formulation to describe two-phase wetting-non-wetting immiscible flow in one dimension. The saturation of each phase is defined as the volumetric fraction of the volume occupied by that phase. Denote the saturation of the wetting phase by $u$, then for a fully saturation porous medium, the saturation of the non-wetting phase is $1 - u$. Let the gravity act in the positive $x$-direction, for each phase, the Darcy-Buckingham law gives
\begin{align}
	v_\alpha & = - \frac{k_{r \alpha} K }{\mu_\alpha} \frac{\partial}{\partial x} (p_\alpha - \rho_\alpha g x) \nonumber \\
	& = - \lambda_\alpha (\frac{\partial p_\alpha }{\partial x} - \rho_\alpha g), \label{eqn:velo}
\end{align}
where $\alpha = { w, n}$ is an index of the wetting and non-wetting phases, $K$ is the intrinsic permeability of the porous medium, $g$ is the gravitational acceleration constant, $k_{r\alpha}, \mu_\alpha, \lambda_\alpha = \frac{k_{r\alpha}}{\mu_\alpha}, p_\alpha, \rho_\alpha$ and $v_\alpha$ are the relative permeability function, viscosity, mobility, pressure, density and volumetric velocity (flux rate across a unit area) of phase $\alpha$, respectively.

Define the total velocity $v_T = v_n + v_w$ and fractional flow rate of the wetting phase $f = \frac{\lambda_w}{\lambda_w + \lambda_n}$, then the velocity of the wetting phase can be expressed by
\begin{align}\label{eqn:vwfrac}
	v_w = f [v_T + \lambda_n  (\frac{\partial}{\partial x} (p_n - p_w) + (\rho_w - \rho_n) g)].
\end{align}
Under non-equilibrium conditions, Stauffer \cite{stauffer1978time}, Hassanizadeh and Gray \cite{hassanizadeh1993thermodynamic}, Kalaydjian \cite{kalaydjian1987macroscopic} proposed that the phases pressure difference $p_n - p_w$ can be written as a function of the equilibrium capillary pressure minus the product of the saturation rate of the wetting phase with a dynamic capillary coefficient $\tau$ [Pa s]:
\begin{align}
	\label{eqn:taumodel}
	p_n - p_w = p_c - \tau \frac{\partial u}{\partial t},
\end{align}
where $p_c$ modeling the capillary pressure under equilibrium condition, is a smooth and decreasing function of saturation $u$, and $\tau$ can be explained as a relaxation time. We refer to \cite{hassanizadeh2002dynamic} for a review of experimental work on dynamic effects in the pressure-saturation relationship.

For the wetting phase the mass conservation equation reads
\begin{align}
\label{eqn:masscons}
\frac{\partial (\phi \rho_w u)}{\partial t} + \frac{\partial}{\partial x} (\rho_w v_w) = \rho_w F_w,
\end{align}
where $\phi$ is the porosity of the porous medium and $F_w$ is  source of wetting phase.

Assuming that $\phi$ and temperature are constant, the phases are incompressible and neglecting the source term, using (\ref{eqn:taumodel}) and substituting (\ref{eqn:vwfrac}) into (\ref{eqn:masscons}) give the MBL equation
\begin{align}
	\label{eqn:fullmbl}
	\frac{\partial u}{\partial t} + \frac{\partial F(u) }{\partial x}= - \frac{\partial}{\partial x} [ H(u) \frac{\partial}{\partial x} (p_c(u) - \tau \frac{\partial u}{\partial t})].
\end{align}
In (\ref{eqn:fullmbl}) the flux $F(u)$ and the capillary induced diffusion \cite{cuesta2006non} $H(u)$ are given by
\begin{align}
	& F(u) = \frac{ 1}{\phi}f(u) [ v_T + \lambda_n(u) (\rho_w - \rho_n) g],\label{eqn:flux} \\
	& H(u) = \frac{1}{\phi} \lambda_n(u) f(u)\label{eqn:diff}.
\end{align}
The fractional flow rate $f(u)$ has a characteristic S-shaped graph. When gravity is included, with different values of $v_T$, the graphs of $F(u)$ are illustrated in Fig. \ref{fig:fluxdiagcomp} (left) and Fig. \ref{fig:fullmodelfh} (left). For a realistic model, the diffusion function $H(u)$ degenerates at $u = 0$ and $1$, see Fig. \ref{fig:fullmodelfh} (right). Since the definitions of relative permeability functions only make sense when $u \in [0, 1]$, in the following we restrict, therefore,  $u \in [0, 1]$.

\subsection{Traveling waves}\label{sec:tw}
Traveling wave (TW) solutions of the MBL equation have been investigated in Refs. \cite{van2007new, spayd2011buckley, van2013travelling}. For the Riemann problem
\begin{equation}
  \label{eqn:riemann}
  u(x, 0) = \left\{ \begin{aligned}
    u_l, \quad & x \leq 0, \\
    u_r, \quad & x > 0,
  \end{aligned}
  \right.
\end{equation}
with different combinations of $(u_l, u_r, \tau)$, the MBL equation may have different types of solutions, for example, the admissible Lax shock, rarefaction wave and undercompressive shock \cite{van2007new, spayd2011buckley, van2013travelling}. In this section, we follow \cite{van2007new} and study the TW solutions of the MBL equation.

To find a traveling wave solution for the MBL equation, we introduce the new variable $\eta = x - st $. Substituting $u(\eta)$ into (\ref{eqn:fullmbl}) results in a third order ordinary differential equation (ODE)
\begin{equation}
\label{eqn:twfull}
 \left\{
    \begin{aligned}
	&-s u' + [F(u)]' = - [H(u) p'_c(u) u']' - s \tau [H(u)u'']', \\
	& u(-\infty) = u_l, \quad u(\infty) = u_r, \quad u_l, u_r \in [0, 1],
  \end{aligned}
 \right.
\end{equation}
where prime denotes differentiation with respect to $\eta$. Integrating this equation over $(\eta, \infty)$ and assuming
\begin{align}
  \label{eqn:odeorder2bc}
 [H(u) (p'_c (u) u' - s \tau  u'')] (\pm\infty) = 0 ,
\end{align}
yields the second-order ODE:
\begin{equation}
	\label{eqn:odeoder2}
\left\{
\begin{aligned}
	& - s(u - u_r) + [F(u) - F(u_r)] = -H(u) p'_c(u) u' - s \tau H(u) u'', \\
   & u(-\infty) = u_l, \quad u(\infty) = u_r,
\end{aligned}
\right.
\end{equation}
with $s$ determined by the Rankine-Hugoniot condition
\begin{align}
	\label{eqn:rhc}
	s = \frac{F(u_l) - F(u_r)}{u_l - u_r}.
\end{align}
When gravity is included into the flux function $F(u)$, Fig. \ref{fig:fluxdiagcomp} (left) and Fig. \ref{fig:fullmodelfh} (left) show that, with different values of $v_T$, $F(u)$ may be non-monotone. For simplicity, we only consider the $(u_l, u_r)$ pairs that satisfy $s > 0$.

In (\ref{eqn:fullmbl}), when $F(u)$, $H(u)$  and $p_c(u)$ are given by
\begin{align}
\label{eqn:symfh}
	F(u) = \frac{u^2}{u^2 + M (1 - u)^2}, \quad H(u) = \epsilon^2, \quad p_c(u) = -\frac{u}{\epsilon},
\end{align}
Van Duijn et al. in Ref. \cite{van2007new} proved that the existence of the TW solution satisfying (\ref{eqn:odeoder2}) depends on the values of $(u_l, u_r, \tau)$.

If we consider (\ref{eqn:symfh}), we can summarize the results obtained by Ref. \cite{van2007new} as follows. Let $u_I$ be the unique inflection point of the flux function $F(u)$. Consider $u_0 \in [0, u_I)$, then it's proved that there is a constant $\tau_*$ such that for all $\tau \in [0, \tau_*]$, there exists a unique solution of (\ref{eqn:odeoder2}) connecting $u_l = u_\alpha$ and $u_r = u_0 $, where $u_\alpha$ is the unique root of the equation
	\begin{align}
		F'(u) = \frac{F(u) - F(u_0)} {u - u_0}.
	\end{align}
When $\tau > \tau_*$, there exists a unique constant $\bar{u} > u_\alpha$, such that (\ref{eqn:odeoder2}) has a unique solution connecting $u_l = \bar{u}$ and $u_r = u_0$. For $u_r = u_0 < u_l = u_B < \bar{u}(\tau)$, the solution of (\ref{eqn:odeoder2}) will exist only if $u_B \in (u_0, \underline{u})$, where $\underline{u}$ is the unique root in the interval $(u_0, \bar{u}$) of
\begin{align}
\label{eqn:uunderline}
	\frac{F(u) - F(u_0)}{u - u_0} = \frac{F(\bar{u}) - F(u_0)} { \bar{u} - u_0}.
\end{align}
When $\tau > \tau_*$ and $u_B \in (\underline{u}, \bar{u})$, there is no TW solution of (\ref{eqn:odeoder2}) connecting $u_l = u_B$ and $u_r = u_0$. In this situation, the solution profile is non-monotonic, two TWs are used in succession: one from $u_l = u_B$ to $u_r = \bar{u}$ and one from $u_l = \bar{u}$ to $u_r = u_0$. For any $u_B \in (\underline{u}, \bar{u})$ and $\tau > \tau_*$, there exists a unique solution of (\ref{eqn:odeoder2}) such that $u_l = u_B$, $u_r = \bar{u}$.

\begin{figure}[!htbp]
\begin{center}
	{\includegraphics[width=2.7in] {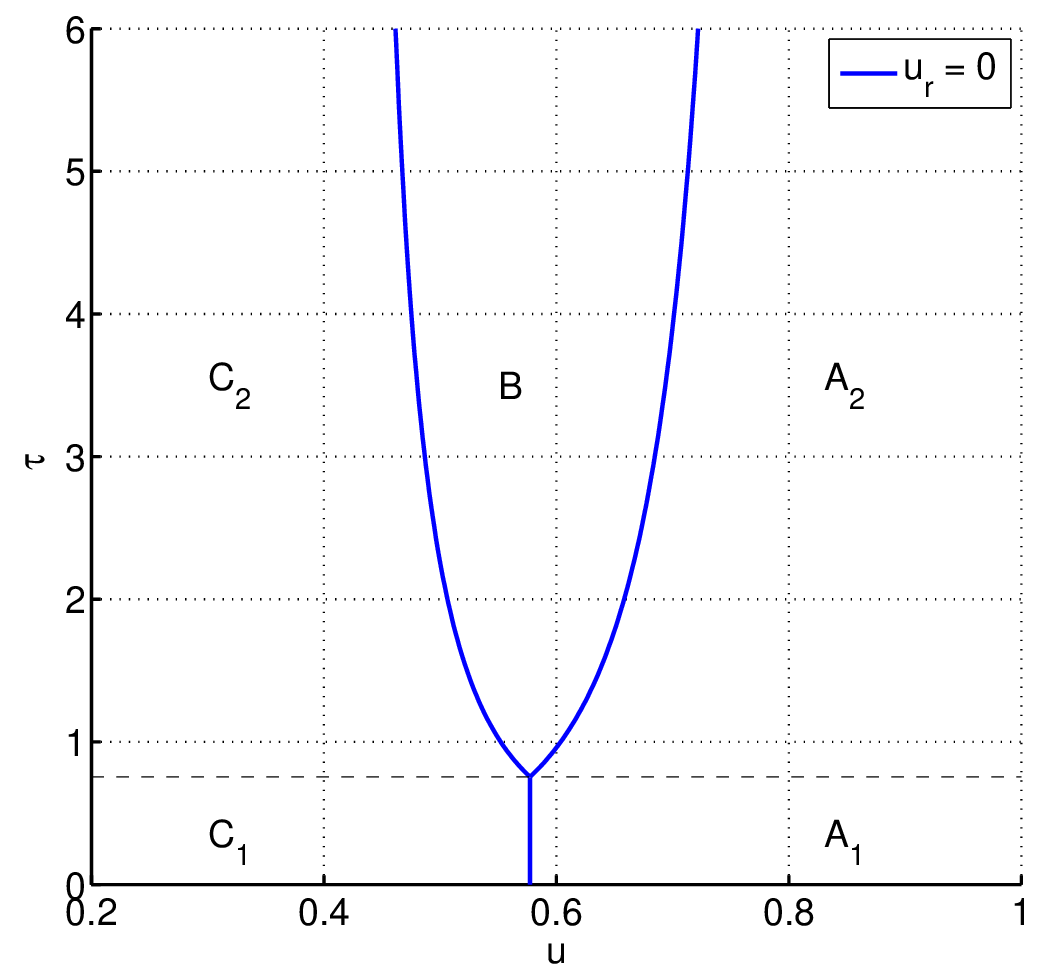}}
	 \caption{Bifurcation diagram for the flux function $f$ with $u_0 = 0$).\label{fig:kaobifurcationregions} }
\end{center}
\end{figure}

For a given $\bar{u} > u_\alpha$, an algorithm to determine the value of $\tau$ is presented in Ref. \cite{van2007new}. This is based on the following concept, invert the function $u(\eta)$ and define the new dependent variable $w(u)= - u'(\eta(u))$, which satisfies
\begin{align}
\label{eqn:slope}
 s \tau H(u)  w w' - H(u) p'_c(u)  w = s ( u - u_r) - [F(u) - F(u_r)],
\end{align}
with boundary condition
\begin{align}
\label{eqn:slopeboundary}
w(u_r = u_0) = w(u_l =  \bar{u}) = 0.
\end{align}
The value of $\tau$ corresponding to a given $\bar{u}$ can be computed using a shooting method. For more details on this, we refer to \cite{van2007new}. To show the relationship between $\tau \text{-} \bar{u}$, we take $M = 0.5$, $\epsilon = 10^{-3}$, and plot the bifurcation diagram for $u_0 = 0$ in Fig. \ref{fig:kaobifurcationregions}.

When $u_0 < u_I$ and $u_B > u_0$, the travelling solutions can be classified using the five regions in the bifurcation diagram. The results summarized from Ref. \cite{van2007new} are given in Table \ref{tab:kaobifurcationregions}.

\begin{table}
  \caption{Results summarized from Ref. \cite{van2007new}\label{tab:kaobifurcationregions}}
  \center
  \begin{tabular}{|c|p{300pt}|}\hline
    Region & Solution description \\ \hline
    $(u_B, \tau) \in A_1$ & Rarefaction wave from $u_B$ down to $u_\alpha$ trailing an admissible Lax shock from $u_\alpha$ down to $u_0$ \\ \hline
    $(u_B, \tau) \in A_2$ & Rarefaction wave from $u_B$ down to $\bar{u}$ trailing an undercompressive shock from $\bar{u}$ down to $u_0$ \\ \hline
    $(u_B, \tau) \in B$ & An admissible Lax shock from $u_B$ up to $\bar{u}$ (may exhibit oscillations near $u_l = u_B$) trailing an undercompressive shock from $\bar{u}$ down to $u_0$ \\ \hline
    $(u_B, \tau) \in C_1$ & An admissible Lax shock from $u_B$ down to $u_0$ \\ \hline
    $(u_B, \tau) \in C_2$ & An admissible Lax shock from $u_B$ down to $u_0$ (may exhibit oscillations near $u_l = u_B$ \\ \hline
  \end{tabular}
\end{table}

Next we write Eq. (\ref{eqn:odeoder2}) as a first order system of ODEs:
\begin{equation}\label{eqn:odesorder1}
	\left\{ \begin{aligned}
		&u' = v, \\
		&v' = \frac{1} {s  \tau H(u)} \big[s(u- u_r) - [F(u) - F(u_r)] - H(u)p'_c(u) v \big].
	\end{aligned}\right.
\end{equation}
When $u_l \neq u_\alpha, u_r = u_0$, the ODE system has three equilibria:
\begin{align}
	(u,v) = (u_0, 0), \quad (u,v) = (\underline{u}, 0), \quad (u, v) = (\bar{u}, 0).
\end{align}
The Jacobian of (\ref{eqn:odesorder1}) reads
\begin{equation}
	A = \left[ \begin{array}{cc}
			0 & 1\\
			\frac{s - F'(u) }{s \tau H(u) } & -\frac{H(u)p_c'(u)}{s \tau  H(u) }
	\end{array} \right],
\end{equation}
and has eigenvalues
\begin{align}
\label{eqn:eigenvalue}
	\lambda_{\pm} = \frac{1}{2 s \tau } [- p'_c(u) \pm \sqrt{(p'_c)^2 - 4 s \tau \frac{ F'(u) -s )}{H(u)}} ].
\end{align}
From this we can get the classification of the three equilibria. The outside two equilibria $(u_0, 0)$ and $(\bar{u}, 0)$ are saddles and the middle equilibrium $(\underline{u},0)$ is either an unstable node or a spiral since $F'(u) > s$.
When $\underline{u} < u_B < \bar{u}$ (or $u_0 < u_B < \underline{u})$, consider a traveling wave connecting $u_l = u_B$ and $u_r = \bar{u}$ (or $u_r = u_0$) and the wave speed is
\begin{align}
  s = \frac{F(u_l) - F(u_r)}{u_l - u_r}.
\end{align}
When using (\ref{eqn:eigenvalue}) and $\tau > \tau_s = \frac{H(u_l) p'_c(u_l)^2 }{4 s (F'(u_l) - s)}$, the equilibrium $(u_B, 0)$ is a spiral.

\textit{Remark} When $u_I < u_0 < 1$ and $0 < u_B < u_0$, the TW solutions of (\ref{eqn:odeoder2}) can be obtained in a similar way. In this case, the roles of $\bar{u}$ and $\underline{u}$ have been switched: $\underline{u}$ is known as the basin height, $\bar{u}$ is the unique root in the interval $(\underline{u}, u_0)$ of 
\begin{align}
	\label{eqn:ubar}
	\frac{F(u) - F(u_0)}{u - u_0} = \frac{F(\underline{u}) - F(u_0)} { \underline{u} - u_0}.
\end{align}
The boundary condition (\ref{eqn:slopeboundary}) is also replaced by
\begin{align}
	w(u_r = u_0) =  w(u_l = \underline{u}) = 0.
\end{align}

\section{The adaptive moving mesh method}
With the appearance of the non-monotonic profiles as mentioned in the previous section, an ideal mesh used in the simulation should be able to capture the overshoot of the saturation on the wetting front. When the initial saturation $u_0$ is very small, the solution of ODE (\ref{eqn:slope}) may have a large magnitude (see Fig. \ref{fig:case4ur003ur03} (right)), which means a sufficiently dense grid should be used near the wetting front to resolve the sharp profile. Note that, the solution is smooth in the region far away from the wetting front, thus a coarse mesh could be used in this region. Based on the above observations, we choose an adaptive moving mesh method \cite{zegeling2007theory,huang1997analysis} to distribute the grid points dynamically in accordance with the evolution of the solution.

Adaptive moving mesh method continuously repositions a fixed number of grid points according to a monitor function $\omega$, so that the resolution in particular locations of the computational domain is improved. Generally, to apply the moving mesh method, three steps have to be taken.
\begin{enumerate}
	\item Transform the PDE from the physical domain to a computational domain.
	\item Define the adaptive mesh transformation.
	\item Discretize the coupled system of PDEs in the spatial direction, then compute the numerical solution by applying a suitable time-integrator to the semi-discrete system.
\end{enumerate}

\subsection{Transformation of the two-phase flow equation}
Let $x$ and $\xi$ denote the physical and computational coordinates. Without lose of generality, $x$ is assumed to be in the interval $\Omega_p = [x_l, x_r]$ and $\xi \in \Omega_c = [0, 1]$.  A general coordinate transformation between $x$ and $\xi$ is given by
\begin{equation}
	x = x(\xi, t), \quad \xi \in [0, 1], t \in [0, T],
\end{equation}
with
\begin{align}
	x(0, t) = x_l, \quad x(1,t) = x_r,	
\end{align}
where $t$ denotes time. In the new coordinate, using the total differential $\dot{u} = \frac{\mathrm{d} u}{\mathrm{d} t} = u_t + u_x \dot{x}$, the physical PDE is transformed to its Lagrangian form
\begin{align}
\label{eqn:transpde}
(\mathcal{I} - \tau \frac{\partial}{\partial x} H(u) \frac{\partial}{\partial x} ) ( \dot{u} - u_x \dot{x}) + \frac{\partial}{\partial x}  F(u) + \frac{\partial}{\partial x} [H(u) \frac{\partial}{\partial x} p_c(u)] = 0.
\end{align}
where $\mathcal{I}$ is the identity operator.
In the next section the transformed PDE will be coupled with a moving mesh PDE (MMPDE) which defines the mesh movement and the monitor function.

\subsection{Mesh transformation with smoothing}
Given a uniform computational mesh with space step $\Delta \xi = \frac{1}{N}$, $\xi_i = \frac{i}{N}, i = 0, 1, \cdots, N$, an adaptive physical mesh $x_i, i = 0, 1, 2, \cdots, N$ is built to equidistribute a specified monitor function $\omega$. In continuous form, the equidistribution principle (EP) \cite{huang1994moving} of the mesh can be expressed as
\begin{align}
\label{eqn:epcontinuous}
	\int_{x_{i-1}}^{x_i} \omega \mathrm{d} x = \int_{x_i}^{x_{i+1}} \omega \mathrm{d} x = c, \quad 1 \leq i \leq N-1,	
\end{align}
or in discrete form
\begin{align}
  \label{eqn:discrete}
  \omega_{i-1} \Delta x_{i-1} = \omega_i \Delta x_i = c, \quad 1 \leq i \leq N-1,
\end{align}
where $\Delta x_i = x_{i+1} - x_i$ is the local grid spacing, $\omega_i$ is a discrete approximation of the monitor function $\omega$ in the grid interval $[x_i, x_{i+1}]$, and $c$ is a constant determined from
\begin{align}
  \label{eqn:constantc}
  \int_{\Omega_p} \omega \mathrm{d} x = \sum_{i = 1}^N \int_{x_{i-1}}^{x_i} \omega \mathrm{d} x = N c.
\end{align}
The monitor function $\omega$ is chosen to cluster mesh points in regions where more accuracy is needed, so it's usually taken to be some measure of the error estimated from the discrete solution. As is often seen in literature, for a scalar solution $u$, a popular choice for controlling grid concentration is based on the arc-length type monitor 
\begin{align}
    \label{eqn:arclengthmonitor}
  \omega = \sqrt{1 + \alpha |u_x|^2},
\end{align}
where the parameter $\alpha$ controls the amount of adaptivity, in this work we set $\alpha = 1$. The choice of the `optimal' monitor function according to interpolation error estimates has been discussed in Ref. \cite{huang2003variational}. In this work, we consider a time-dependent monitor function proposed by Ref. \cite{zegeling2005robust}
\begin{align}
  \label{eqn:monitor}
  \omega = (1 - \beta) \alpha(t) + \beta |u_\xi|^{\frac{1}{m}},
\end{align}
where the intensity controlling parameter $\alpha(t)$ is defined as
\begin{align}
  \alpha(t) = \int_{\Omega_c} |u_\xi|^{\frac{1}{m}} \mathrm{d} \xi.
\end{align}
In (\ref{eqn:monitor}), the critical regions are identified by the computational derivative $u_\xi$, which is smoother than the physical derivative $u_x$. The adaptivity function $\alpha(t) > 0$ averages the derivative $u_\xi$, resulting in a time-dependent monitor function. In this work, we take $m = 1$, this choice is verified to be robust and efficient in Refs. \cite{zegeling2005robust, van2010balanced}.

Following the approach from \cite{huang2001practical}, we can derive that
\begin{align}
   \int_{\Omega_c} \omega \mathrm{d}\xi  &= \int_{\Omega_c} [(1 - \beta) \int_{\Omega_c} |u_\xi|^\frac{1}{m} \mathrm{d} \xi + \beta |u_\xi|^\frac{1}{m} ] \mathrm{d} \xi \\
	 & = \int_{\Omega_c} [(1 - \beta) |u_\xi|^\frac{1}{m} + \beta|u_\xi|^\frac{1}{m}] \mathrm{d} \xi \\
  & = \int_{\Omega_c} |u_\xi|^\frac{1}{m} \mathrm{d} \xi,
\end{align}
thus $\beta$ is indeed the ratio of points in the critical areas. In this paper we choose $\beta = 0.9$, which means approximately $90\%$ grid points are distributed in the critical regions \cite{huang2001practical}. 

The accuracy of the spatial derivative approximations and stiffness of the system after the space discretization are largely influenced by the regularity of the mesh. To equidistribute the monitor function, we adopt a MMPDE with smoothing \cite{zegeling2007theory,huang1997analysis},
\begin{equation}\label{eqn:smmmpde}
\left\{
\begin{aligned}
	& \frac{\partial}{\partial \xi} \left( \frac{\dot {\tilde{n}}}{\omega}\right) = -\frac{1}{\tau_s} \frac{\partial }{\partial \xi} \left( \frac{\tilde{n}}{\omega}\right), \\
	&\tilde{n} = [\mathcal{I} - \sigma_s (\sigma_s + 1) (\Delta \xi)^2 \frac{\partial ^2}{\partial \xi^2}] n,
\end{aligned}
\right.
\end{equation}
where $\sigma_s$ and $\tau_s$ are the spatial and temporal smoothing parameters, $n = \frac{1}{x_\xi}$ is the point concentration, $\Delta \xi$ is  the space step of the computational domain after discretization.

Refs. \cite{zegeling2007theory, huang1997analysis} show that this smoothed MMPDE has the following properties
\begin{enumerate}
	\item No node-crossing: $J = x_\xi > 0$, in discrete version it reads, $\Delta x_i(t) > 0, \quad \forall t \in [0, T]$.
	\item Local quasi-uniformity: $|\frac{x_{\xi \xi}}{x_\xi}| \leq \frac{1}{\sqrt{\sigma(\sigma+1)\Delta \xi}}$ with discretized version:
	\begin{align}
		\frac{\sigma}{\sigma+1} \leq \frac{\Delta x_{i+1}(t)}{\Delta x_i(t)} \leq \frac{\sigma+1}{\sigma}, \quad \forall t \in [0, T].
	\end{align}
	\item When $\sigma_s = \tau_s = 0$ (no smoothing), (\ref{eqn:smmmpde}) fulfills the basic equidistribution principle of the monitor function:
 \begin{align}
   \omega x_{\xi} = \mathrm{constant}, \forall t \in[0, T],
 \end{align}
 in discretized form it reads
	\begin{align}
		\omega_i \Delta x_i  = \text{constant}, \quad \forall t \in [0, T].
	\end{align}
\end{enumerate}
For the choice of the parameters $\tau_s$ and $\sigma_s$, we follow the suggestions in Ref. \cite{zegeling2007theory}. In practice, the choice of the temporal smoothing parameter depends on the timescale in the model: $\tau_s = 10^{-3} \times$ `timescale in PDE model'. The spatial smoothing parameter $\sigma_s$ can be taken as $\sigma_s = \mathcal{O}(1)$. In Section \ref{sec:numerical}, for all numerical experiments using the moving mesh method, we set $\tau_s = 10^{-3}\times T_{end}$ and $\sigma_s = 2$.

\subsection{Discretization of the coupled PDEs}
We employ a finite difference method to discretize the coupled system. Applying the second order centered finite difference scheme in space direction to (\ref{eqn:smmmpde}) yields
\begin{equation}
\left\{ \begin{aligned}
	&\frac{[\mathcal{I} - \sigma_s(\sigma_s+1) \delta_{xx}] (\dot{x}_{i+1} - \dot{x_i})} {\omega_{i+1/2}(x_{i+1} - x_i)^2} - \frac{[\mathcal{I} - \sigma_s(\sigma_s+1)\delta_{xx}] (\dot{x}_{i} - \dot{x}_{i-1})} {\omega_{i - 1/2}(x_i - x_{i-1})^2} = \\
	&\hspace{3cm}   \frac{1}{\tau_s} \left[ \frac{[\mathcal{I} - \sigma_s(\sigma_s+1)\delta_{xx}]\frac{1}{x_{i+1} - x_i}}{\omega_{i+1/2}}- \frac{ [\mathcal{I} - \sigma_s(\sigma_s + 1)\delta_{xx}]\frac{1}{x_i - x_{i-1}}}{\omega_{i - 1/2}} \right], \quad i = 2, 3, \cdots, N-2,\\
	&\dot{x}_{i+1} - 2\dot{x}_i + \dot{x}_{i-1} = 0, \quad i = 1, N-1,\\
	& \dot{x}_ 0 = \dot{x}_N = 0,
\end{aligned}
\right.
\end{equation}
where $\delta_{xx}$ is the second-order difference operator and $\omega_{i+1/2} = (1 - \beta)\alpha(t) + \beta |\frac{u_{i+1} - u_i}{\Delta \xi}|$. The derivative of the point concentration appears in (\ref{eqn:smmmpde}) is discretized as
\begin{align}
  \dot{n}_i = -\frac{\dot{x}_{i+1} - \dot{x_i}}{(x_{i+1} - x_i)^2}, \quad i = 0, 1, \cdots, N-1.
\end{align}

The transformed physical PDE (\ref{eqn:transpde}) is discretized in the same way. Following the method-of-lines approach, the time-integration of the resulting coupled semi-discretized system is solved using the BDF integrator ode15i of Matlab \cite{matlabr2014a}. 

\section{Numerical experiments}\label{sec:numerical}
In this section, we present some numerical results obtained by the moving mesh method presented in the previous section. The first three examples show the accuracy and features of the moving mesh method, the fourth and fifth examples study the effects of the flux rate and initial saturation by taking gravity into account, the last one solves the full equation by utilizing the Brooks-Corey model \cite{brooks1966properties}. 

\subsection{The accuracy of the moving mesh method}
First, we use three examples to test the accuracy of the moving mesh method. Examples 1, 2 and 3 are modifications of the test cases in Ref. \cite{kao2015fast}.

Consider $F(u), H(u)$ and $p_c(u)$ given in (\ref{eqn:symfh}) with $M = 0.5, \epsilon = 10^{-3}$ and initial condition
\begin{equation}
	u(x,0) = \left\{ \begin{array}{ll}
		u_1, \quad & x \in [0, 0.75], \\
		u_2,\quad & x\in (0.75, 2.25), \\
		0, \quad  & x \in [2.25, 3].
	\end{array} 	
	\right.
\end{equation}
With different combinations of $(u_B, u_0, \tau)$, the TW results obtained from Section \ref{sec:tw} are shown in Table \ref{tab:basinplateau}. For this problem, we set the final time $T = 0.5$.

\textbf{Example 1.} $\tau = 3.5, u_1 = 0.25, u_2 = 0.85 $.

In the left part of the initial condition, we have $u_B = u_1$, $u_0 = u_2$, Table \ref{tab:basinplateau} shows $\tau > \tau_s > \tau_*$ and $\underline{u} < u_B = u_1 < \bar{u}$, thus the left part of the solution consists of an admissible Lax shock from $u_B$ down to $\underline{u}$ (with oscillations near $u_B$) and an undercompressive shock from $\underline{u}$ up to $u_0$. In the right part, $u_0 = 0$, since $u_B = u_2 > \bar{u}$ and  $\tau > \tau_*$, the right part of the solution consists of a rarefaction wave from $u_B$ down to $\bar{u}$ and an undercompressive shock from $\bar{u}$ down to $u_0$. In the following, the left and right parts of the solution are called as non-monotone basin and monotone plateau, respectively. The basin height is $\underline{u} = 0.1036$ and the plateau height is $\bar{u}_p = 0.6938$.

\begin{figure}[!htbp]
\begin{center}
	{\includegraphics[width=2.7in] {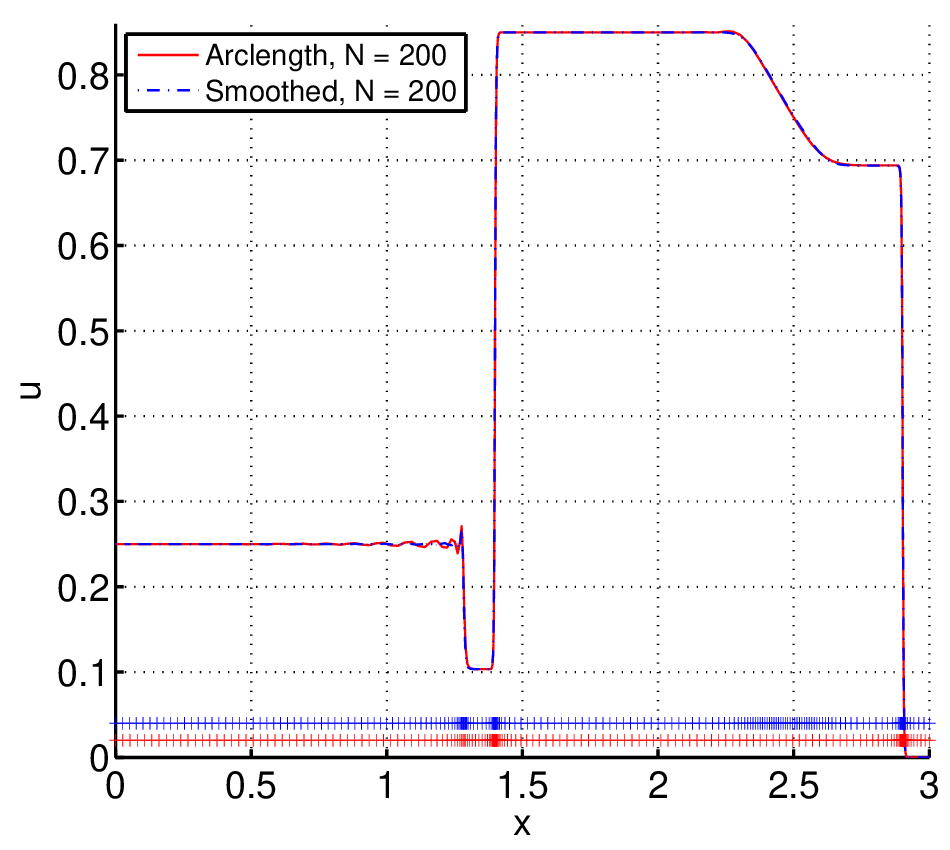}}
    \quad \quad
 {\includegraphics[width=2.7in] {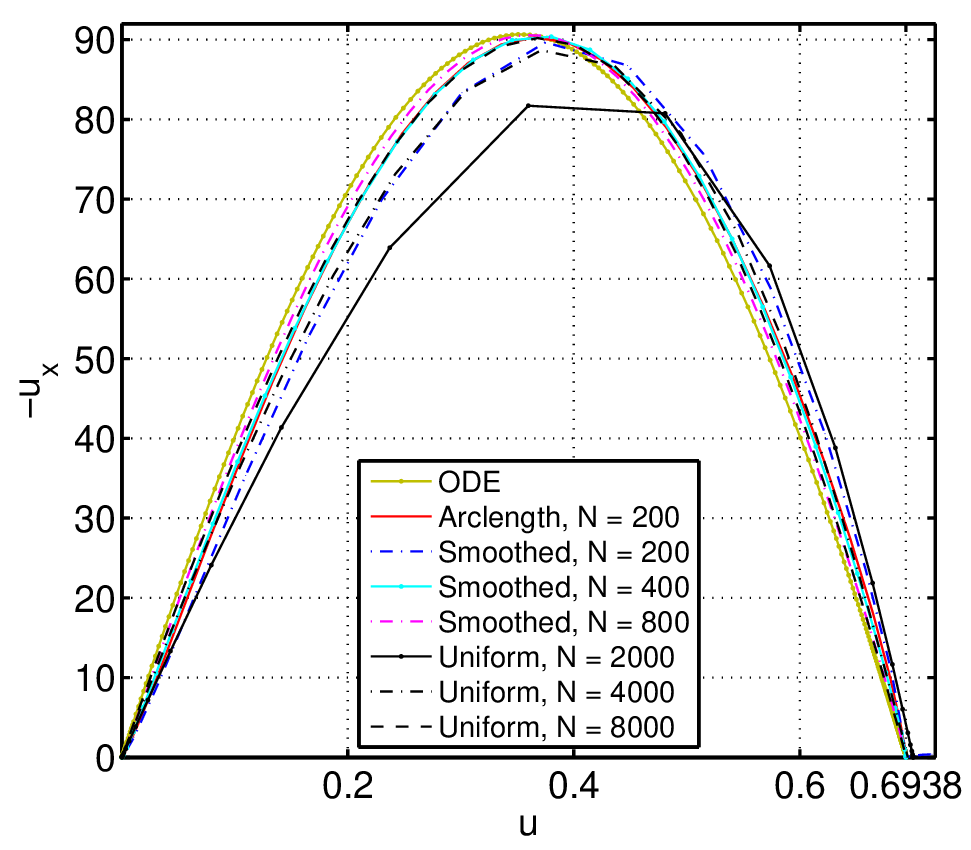}}
    {\includegraphics[width=2.7in] {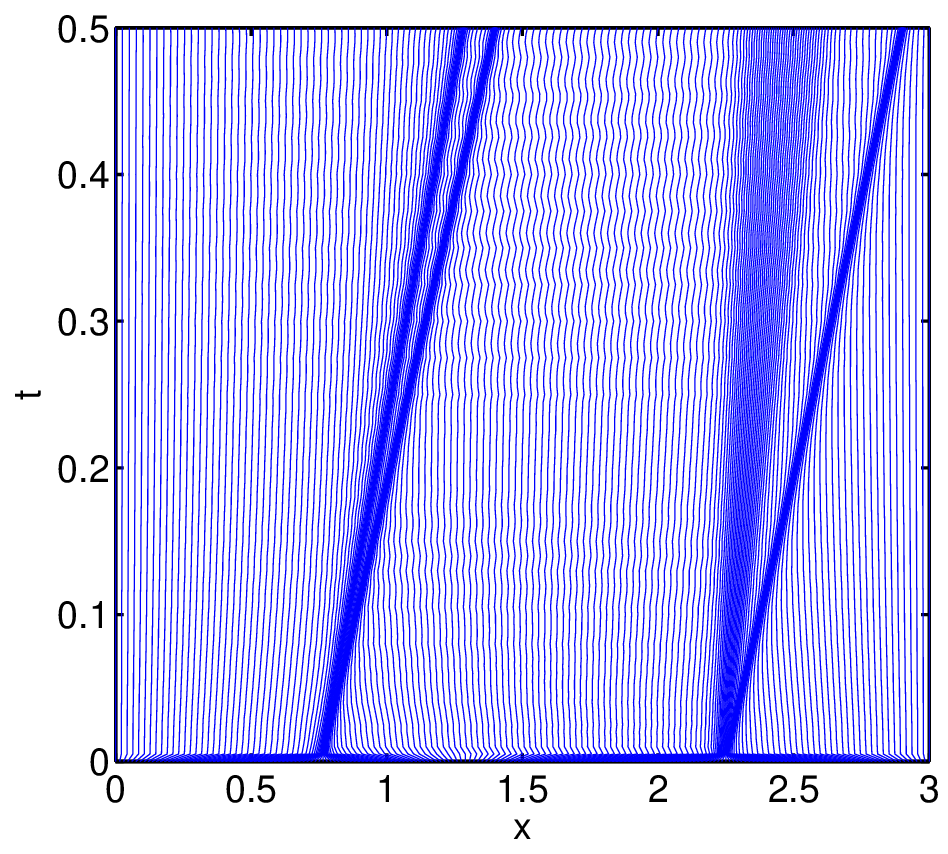}}
    \quad \quad
    {\includegraphics[width=2.7in] {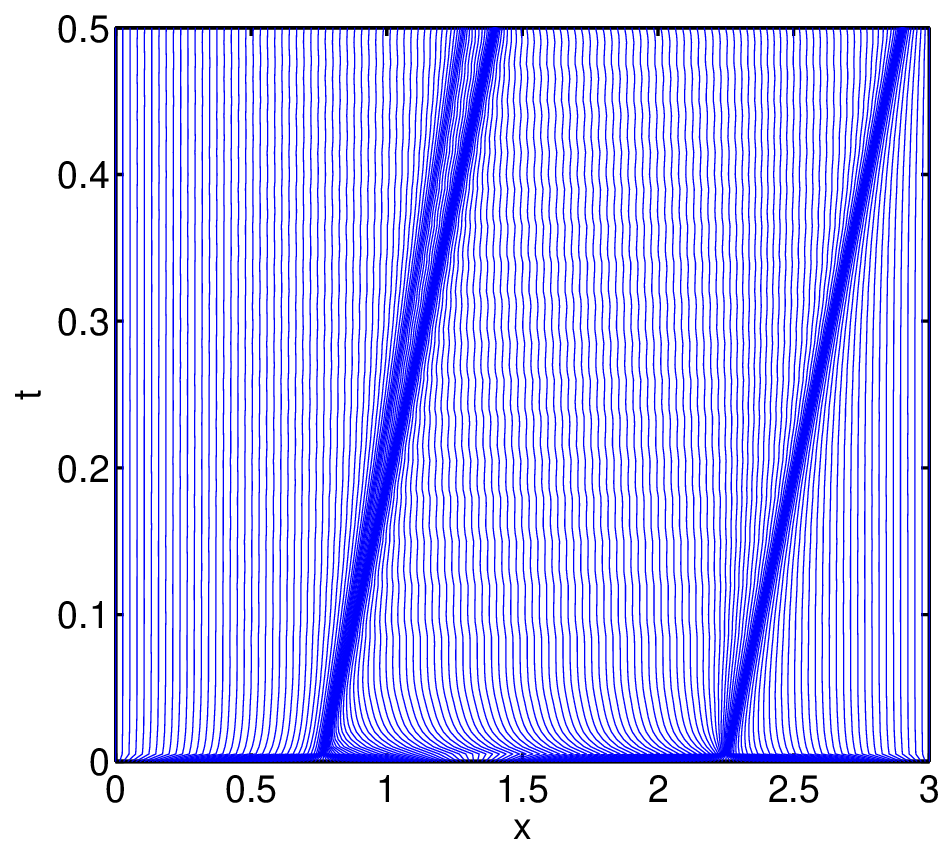}}
   {\includegraphics[width=2.7in] {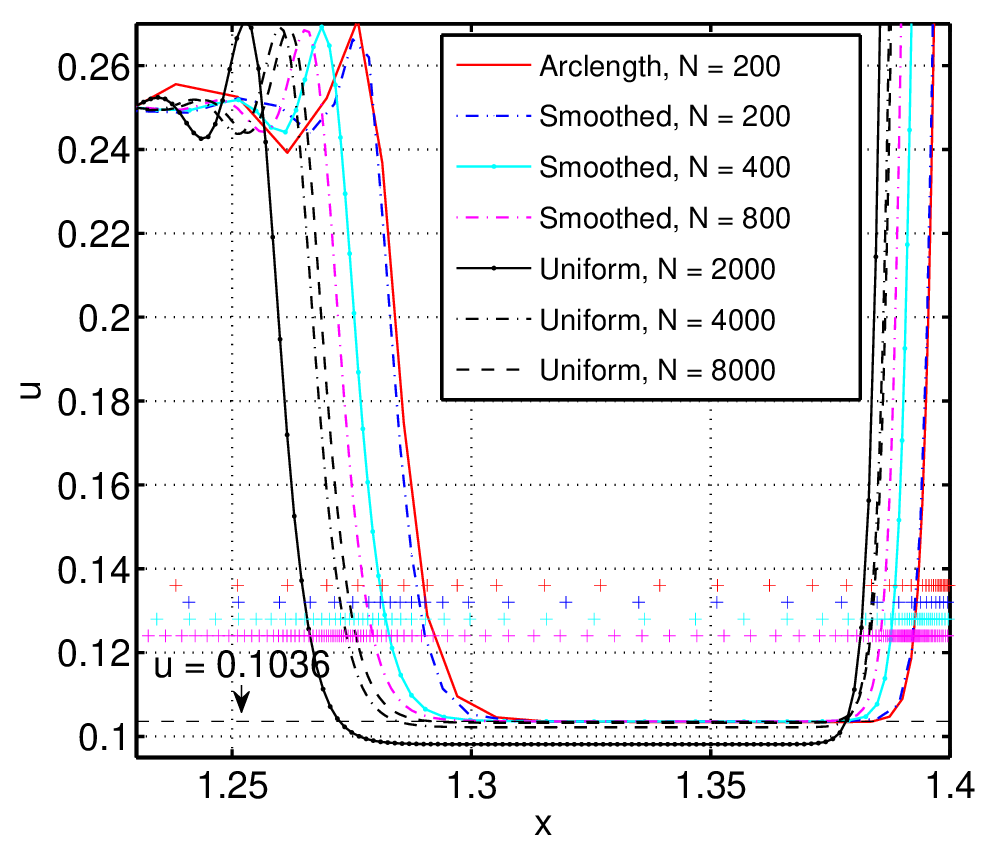}}
   \quad \quad
	 {\includegraphics[width=2.7in] {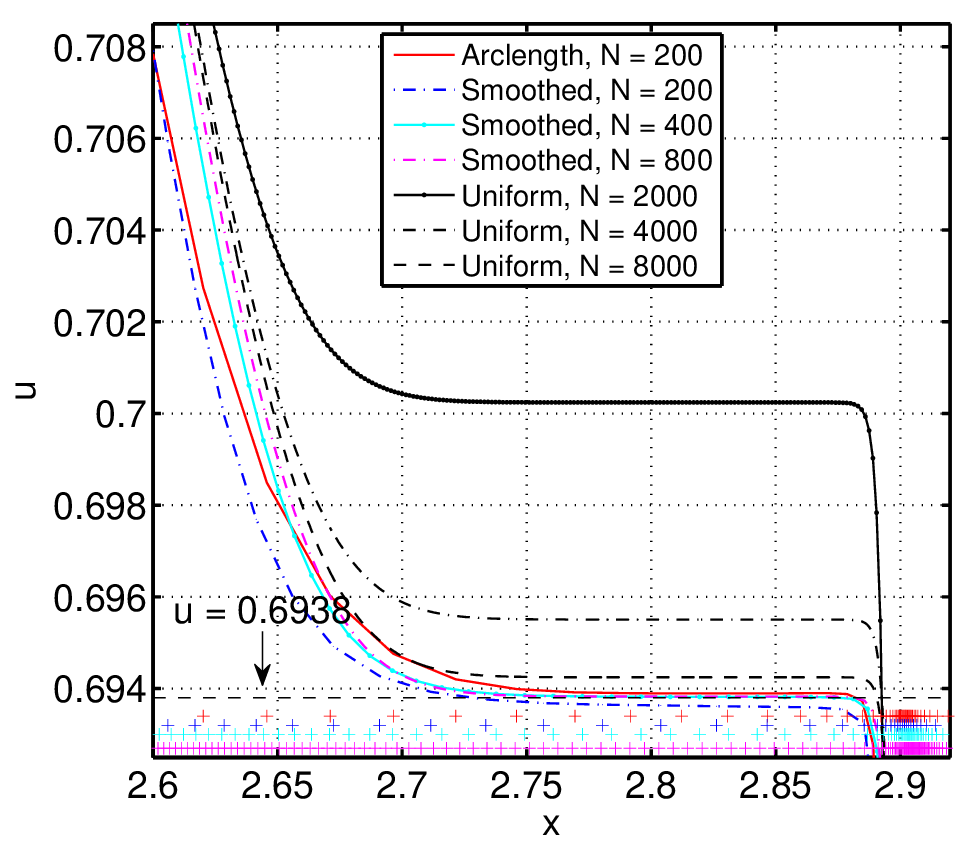}}
	 \caption{Example 1: solutions computed using the moving mesh method (top left); values of $-u_x$ at the right boundary of the plateau (top right); grid trajectories using the smoothed monitor function (middle left) and the arc-length monitor function (middle right) with $N = 200$; zoom in at the basin area (bottom left); zoom in at the plateau area (bottom right). \label{fig:kaoex1}}
\end{center}
\end{figure}

Fig. \ref{fig:kaoex1} shows the computed results obtained by the moving mesh method (monitor (\ref{eqn:arclengthmonitor}) with $N = 200$, monitor (\ref{eqn:monitor}) with $N = 200, 400, 800$) and uniform grids ($N = 2000, 4000, 8000$). Fig. \ref{fig:kaoex1} (top left) clearly shows that the solution includes a non-monotone basin in the left part and a monotone plateau in the right part. Fig. \ref{fig:kaoex1} (top right) plots the opposite slopes ($-u_x$) at the boundary of the right undercompressive shock. The reference slopes are obtained by solving ODE (\ref{eqn:slope}) with the built-in function ode15i in Mablab.  The slopes computed by the moving mesh method using the smoothed monitor (\ref{eqn:monitor}) with $N = 200, 400, 800$ are more accurate than those obtained by the uniform grids with $N = 2000, 4000, 8000$, respectively. The grid trajectories produced by the smoothed monitor (\ref{eqn:monitor}) clearly illustrate the evolution of the solutions. As can be seen, the smoothed monitor attracts more grid points in the rarefaction fan than the arc-length monitor, while the arc-length monitor attracts more grid points near the steep shocks. Therefore, the slope computed using the arc-length monitor is slightly more accurate than the one using the smoothed monitor.

In order to check the accuracy of the moving mesh method, we present the details of the critical areas (basin and plateau) in Fig. \ref{fig:kaoex1} (bottom). The heights of the basin and plateau computed with the moving mesh method ($N = 200, 400, 800$) are more accurate than the heights obtained using a uniform grid with $N = 2000, 4000$ or even $ 8000$. Although the arc-length monitor results in more accurate basin and plateau heights than the smoothed monitor, the accuracy near the admissible Lax shocks is lower as a result of fewer grid points near smooth parts.

\textbf{Example 2.} $\tau = 5, u_1 = 0.25, u_2 = 0.66$.

For the second example, $u_1$ is the same as that in Example 1, $u_2$ is decreased form $0.85$ to $0.66$ and $\tau$ is increased to $5$. Table \ref{tab:basinplateau} shows that this combination results in a non-monotone basin of height $\underline{u} = 0.2027$ in the left part and a non-monotone plateau of height $\bar{u} = 0.7130$ in the right part. In the left part, as a result of the decrease in $u_2$, the height of the basin is higher than the one in Example 1 and the amplitude of the oscillations becomes smaller. In the right part, since $\underline{u} < u_B = u_2 < \bar{u}$ and $\tau > \tau_s > \tau_*$, we obtain a non-monotone plateau which is higher than the one in Example 1. Because $u_B$ is a spiral, a slight oscillation appears near $u_B$.

\begin{figure}[!htbp]
\begin{center}
   {\includegraphics[width=2.7in] {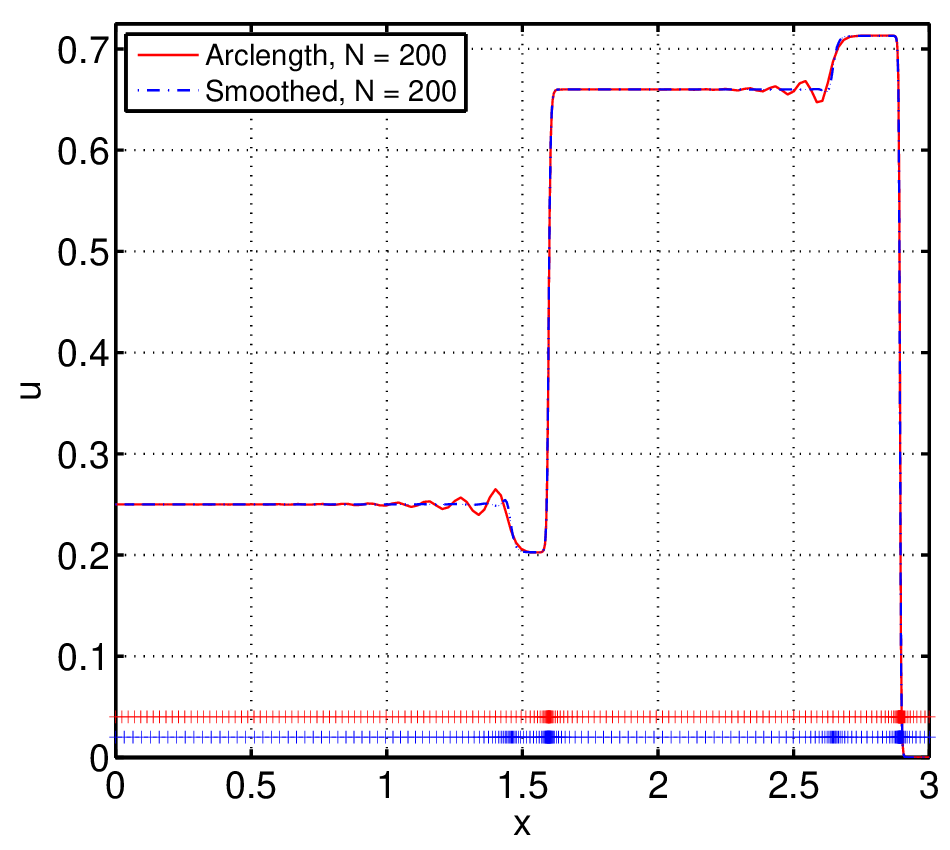}} \quad \quad
    {\includegraphics[width=2.7in] {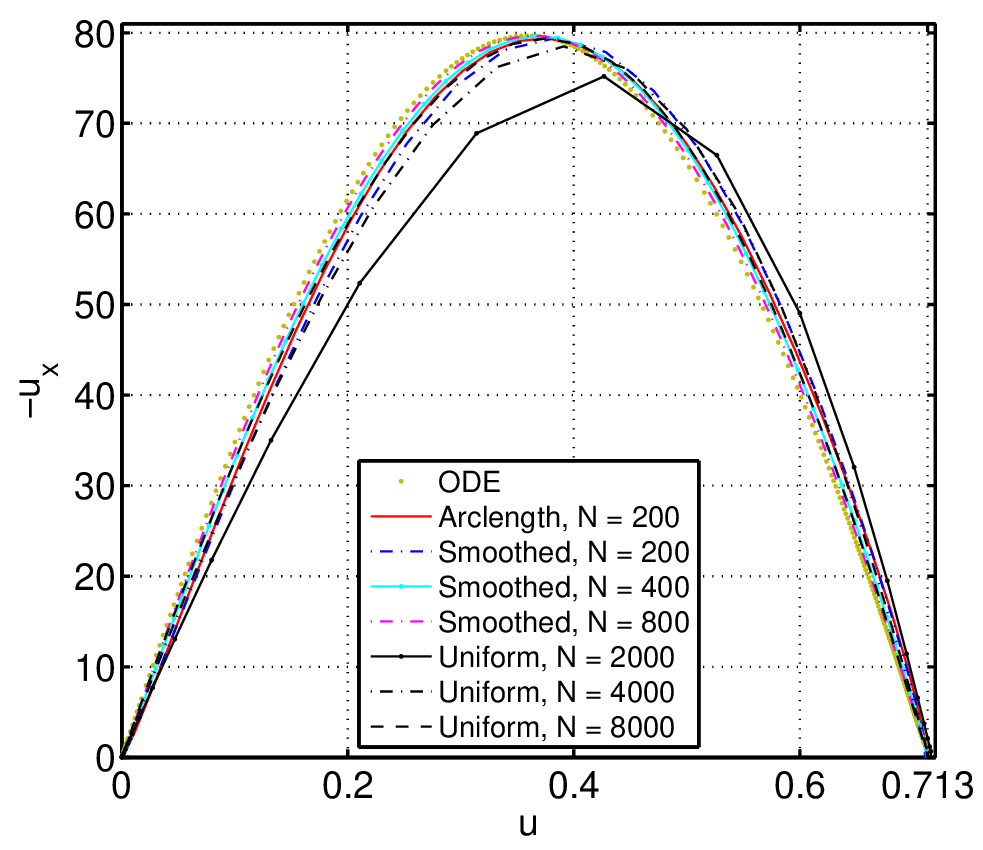}}
   {\includegraphics[width=2.7in] {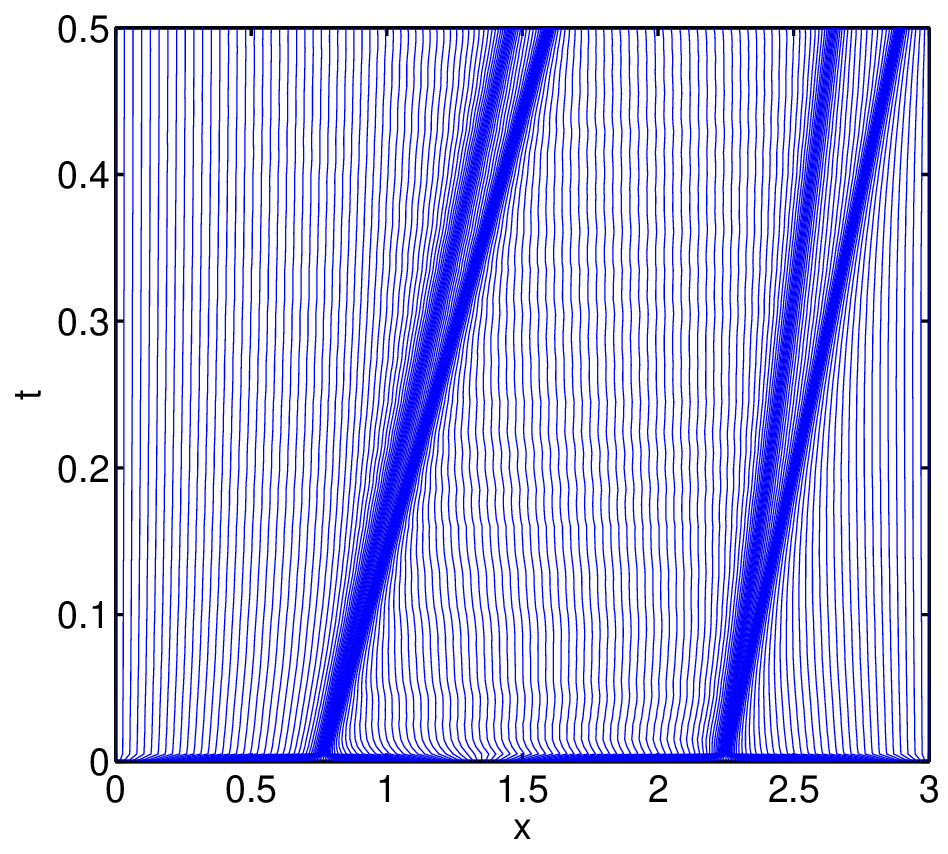}}
   \quad \quad
   {\includegraphics[width=2.7in] {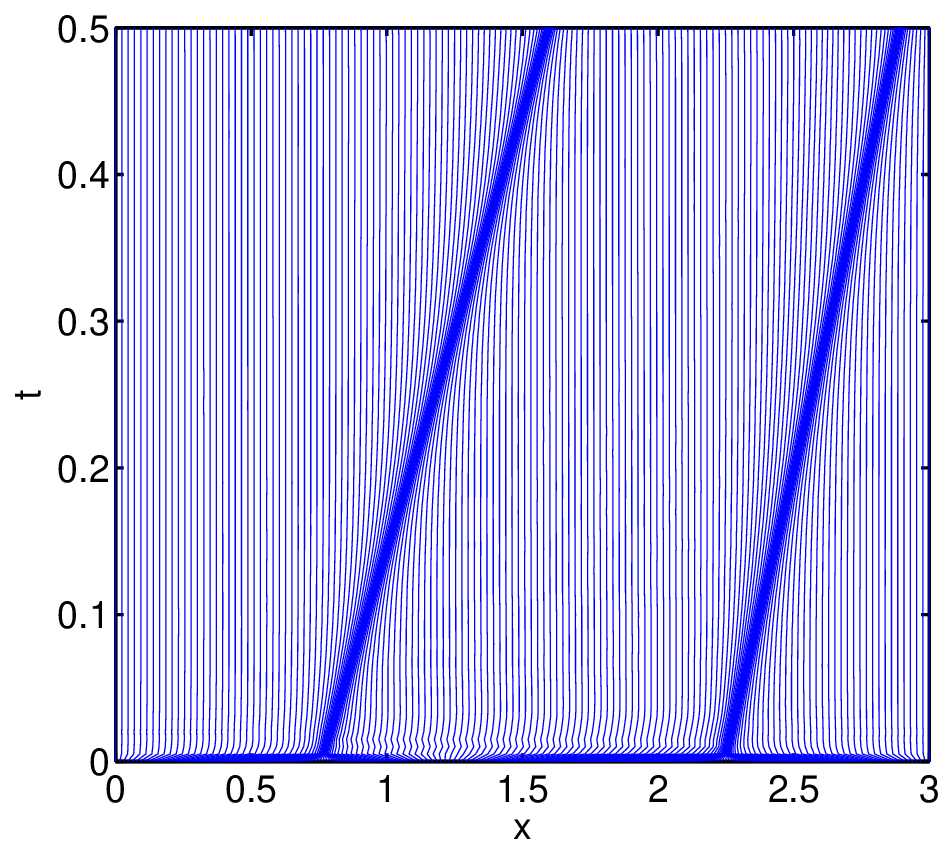}}
	{\includegraphics[width=2.7in] {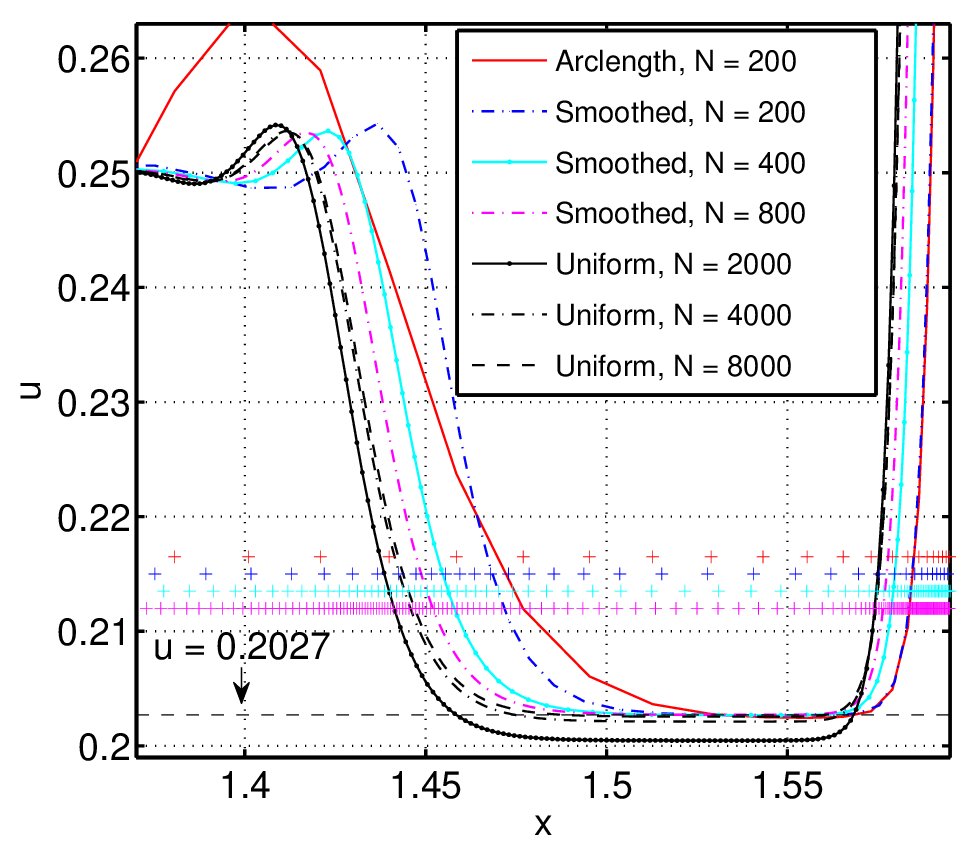}}
    \quad \quad
   {\includegraphics[width=2.7in] {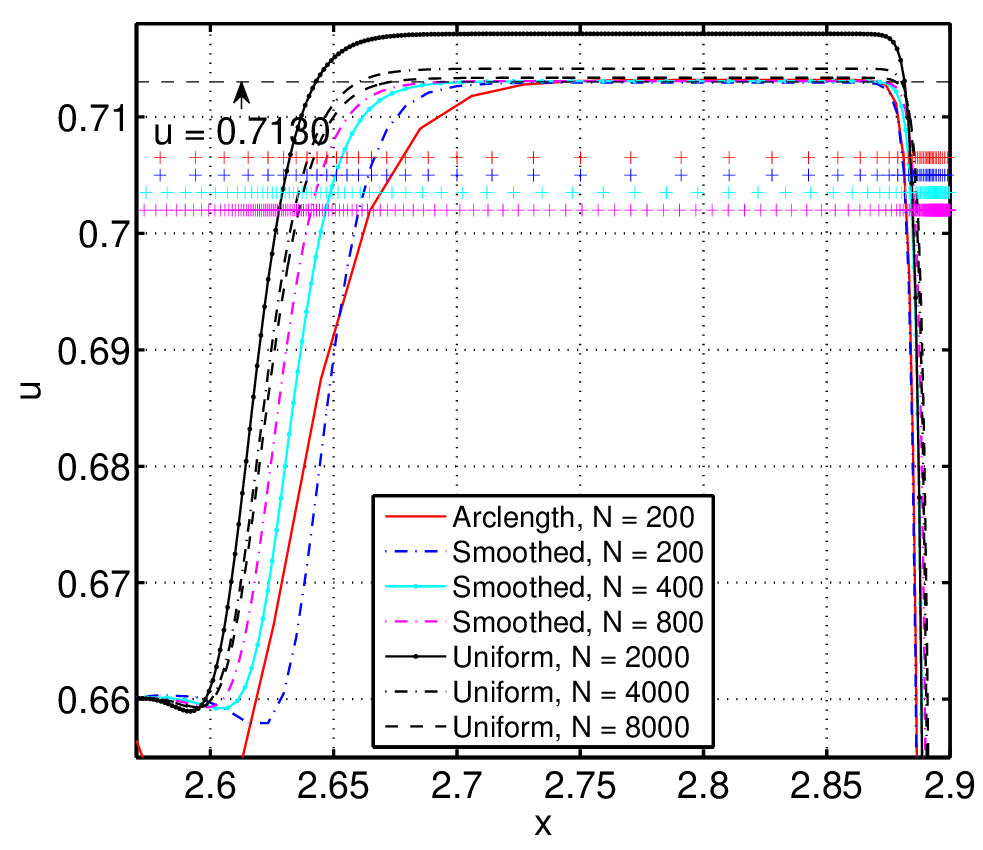}}
	 \caption{Example 2: solutions computed using the moving mesh method (top left); values of $-u_x$ at the right boundary of the plateau (top right); grid trajectories using the smoothed monitor function (middle left) and the arc-length monitor function (middle right) with $N = 200$; zoom in at the basin area (bottom left); zoom in at the plateau area (bottom right). \label{fig:kaoex2}}
\end{center}
\end{figure}

Once again, near the admissible Lax shocks, the smoothed monitor performs better than the arc-length monitor as it attracts more grid points near these parts. The moving mesh method leads to much more accurate slopes, plateau and basin heights than the uniform grid method, see Fig. \ref{fig:kaoex2}.

\textbf{Example 3.} $\tau = 5,  u_1 = 0.25, u_2 = 0.52$.

\begin{figure}[!htbp]
\begin{center}
   {\includegraphics[width=2.7in] {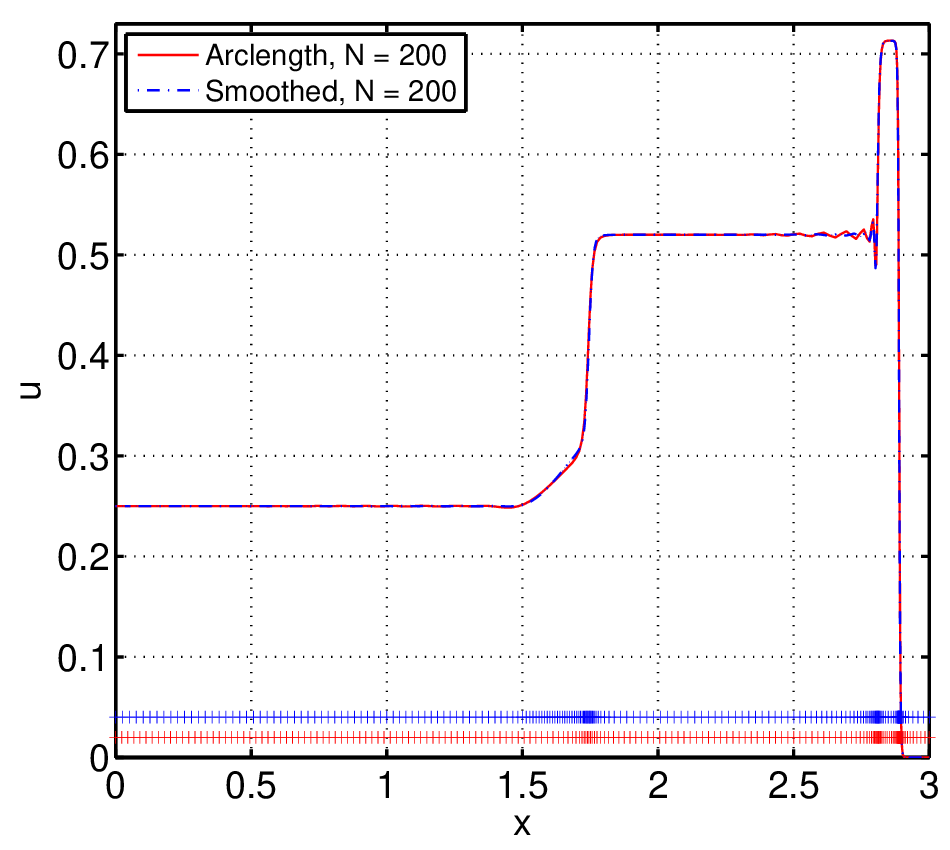}}
   \quad \quad
   {\includegraphics[width=2.7in] {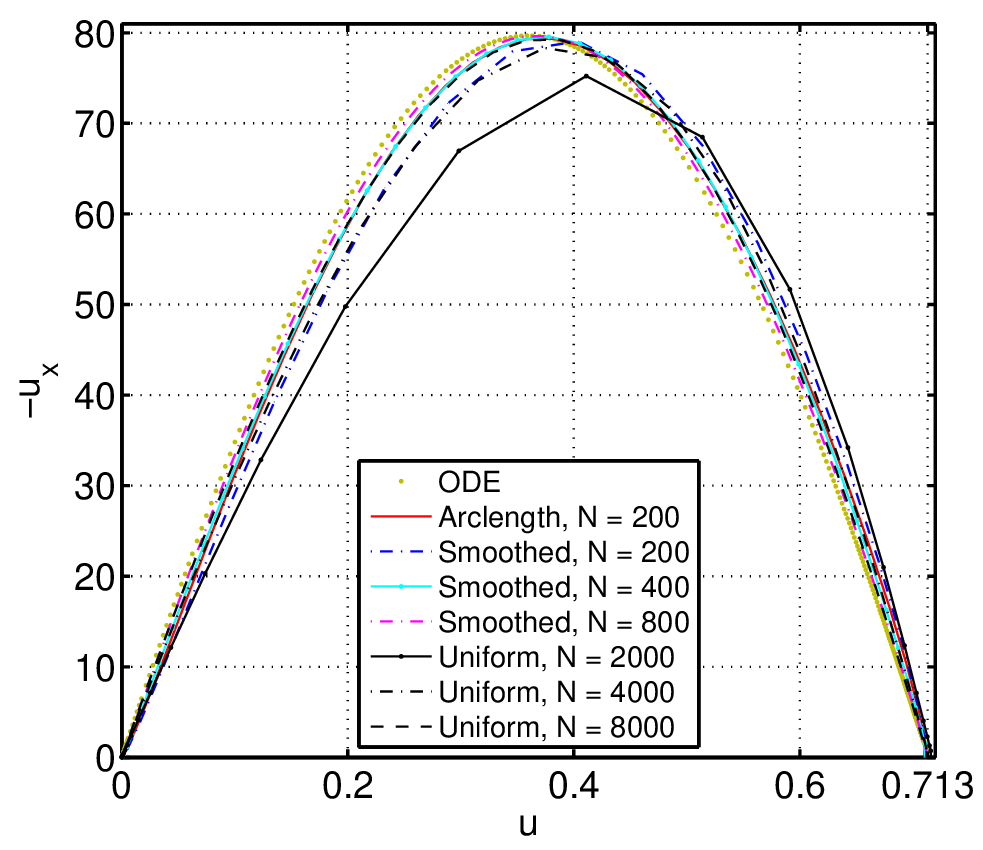}}
   {\includegraphics[width=2.7in] {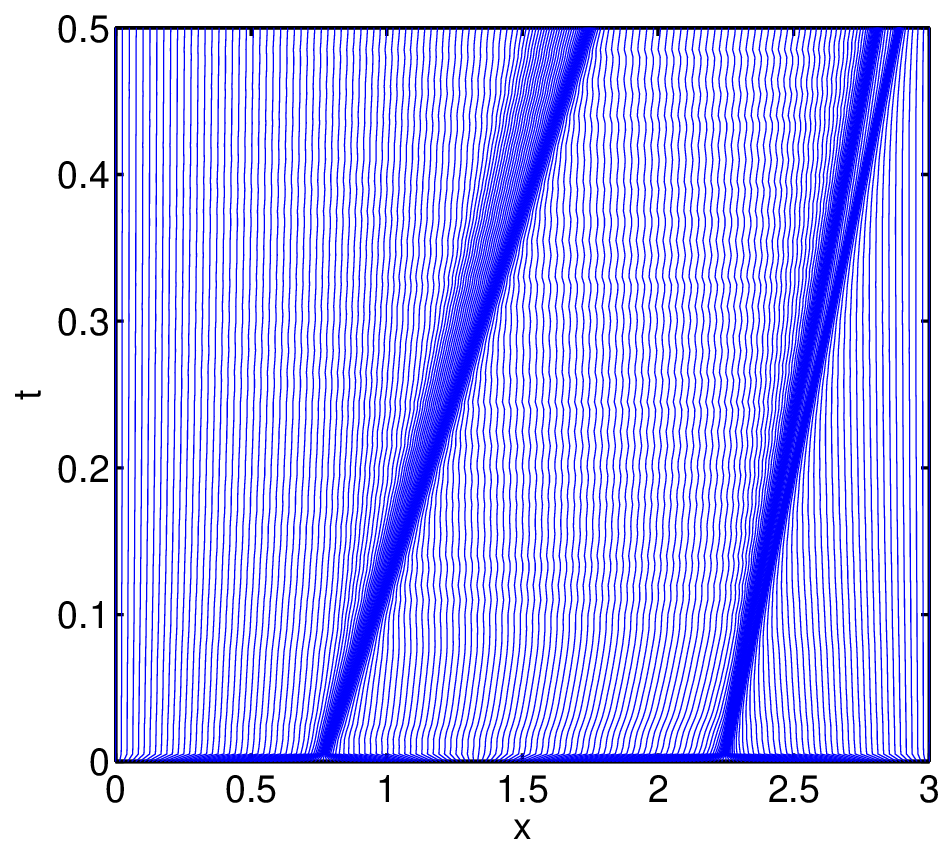}}
    \quad \quad
   {\includegraphics[width=2.7in] {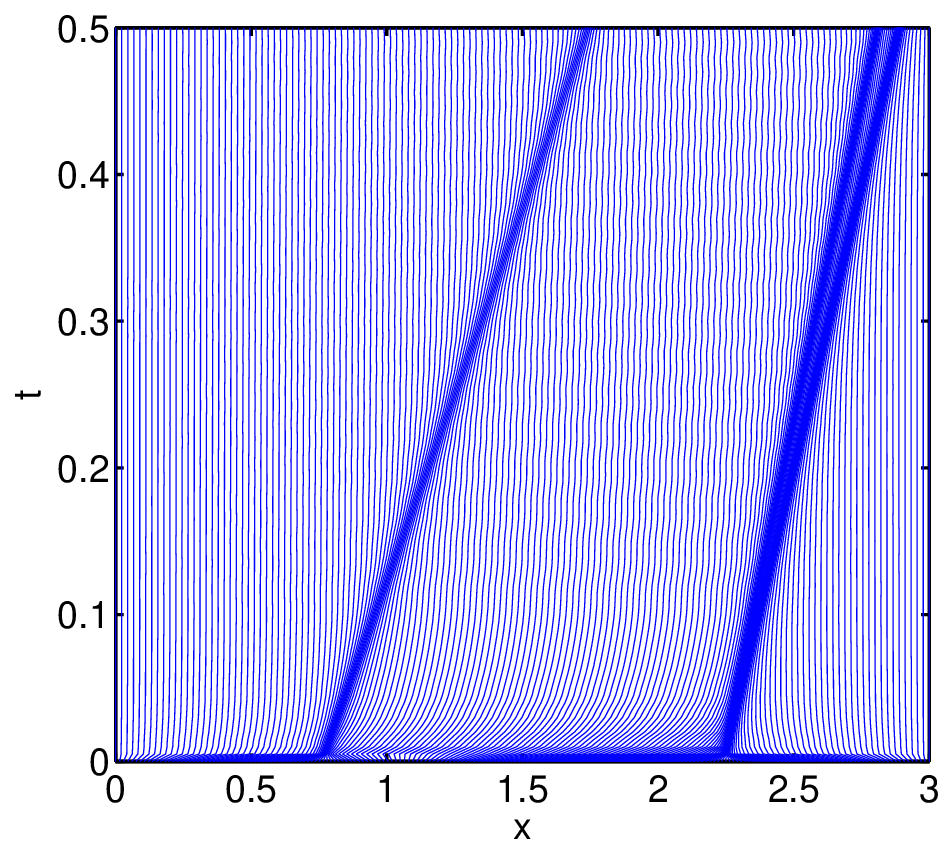}}
 {\includegraphics[width=2.7in] {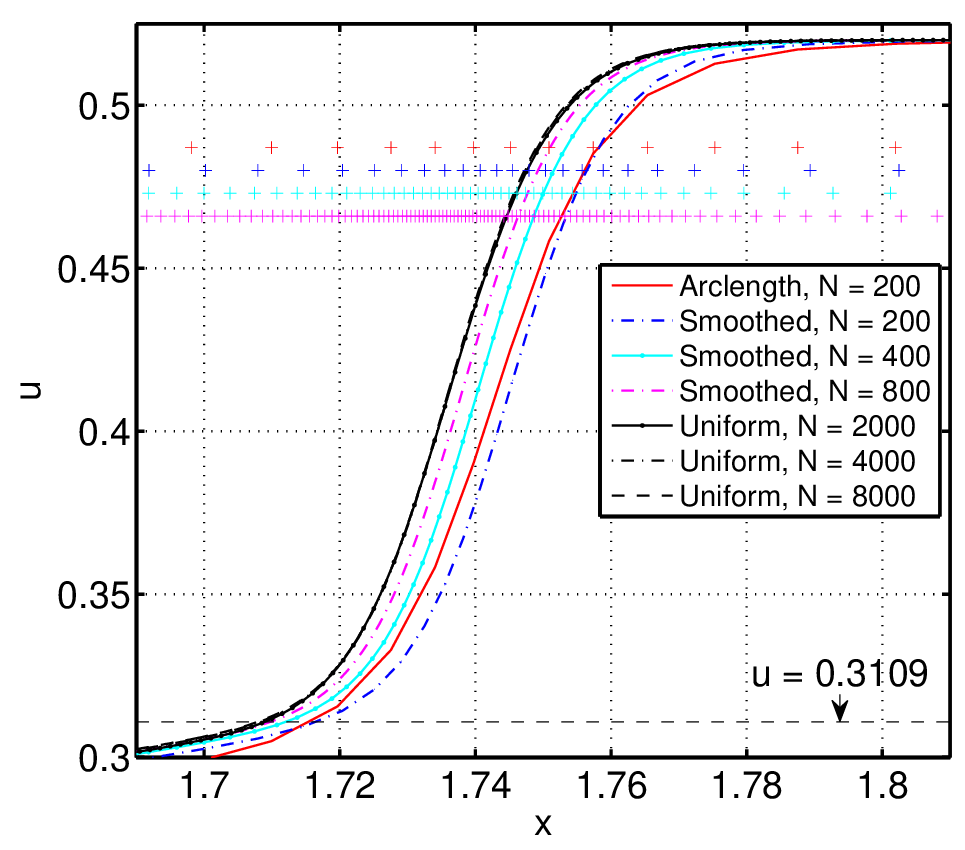}}
    \quad \quad
   {\includegraphics[width=2.7in] {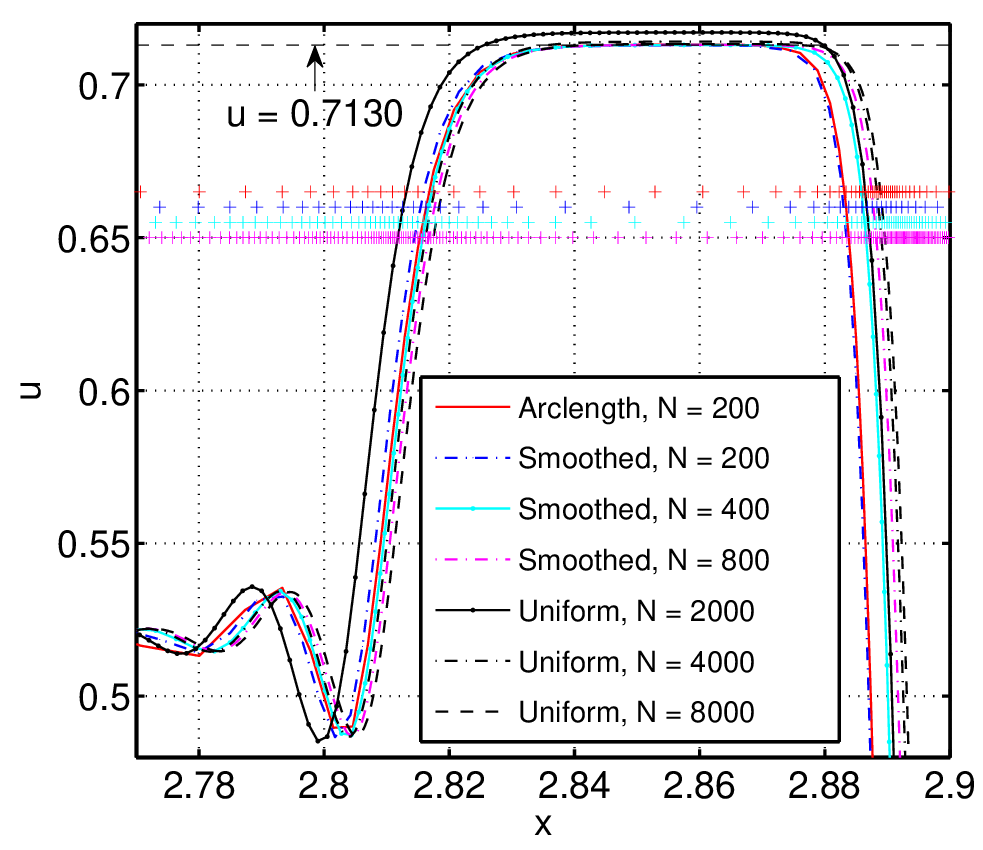}}
	 \caption{Example 3: solutions computed using the moving mesh method (top left); values of $-u_x$ at the right boundary of the plateau (top right); grid trajectories using the smoothed monitor function (middle left) and the arc-length monitor function (middle right) with $N = 200$; zoom in at the basin area (bottom left); zoom in at the plateau area (bottom right). \label{fig:kaoex3}}
\end{center}
\end{figure}

In the third example, a smaller value of $u_2$ is used. In the left part, $u_B = u_1 < \underline{u}$, therefore the solution has a monotone basin area. Note that, since we only computed to $t = 0.5$, instead of a horizonal basin, only a turning point appears near $\underline{u} = 0.3109$. In the right part, $\underline{u} < u_B < \bar{u}$, $\tau > \tau_* > \tau_s$, hence $(u_B, 0)$ is a spiral of the ODE system (\ref{eqn:odesorder1}), consequently oscillations appear near $u_B = u_2$. As in Examples 1 and 2, the plateau heights and slopes obtained by the moving mesh method are more accurate than the uniform grid method, see Fig. \ref{fig:kaoex3}.

\begin{table}
  \caption{\label{tab:basinplateau}Travelling wave results for Example 1, 2 and 3}
  \center
  \begin{tabular}{|c|cccccccc|c|}\hline
		& $u_B$ & $u_0$ & $\tau$ & $\tau_*$ & $\tau_s$ & $u_*$  & $\underline{u}$ & $\bar{u}$ & Wave description \\ \hline
   & 		0.25 & 0.85 & 3.5		 & 0.6826	 & 0.4495 & 0.2151 & 0.1036 & 0.3155 & Non-monotone basin \\
\raisebox{1.6ex}[0pt]{Example 1}& 0.85 & 0 & 3.5 & 0.7545 & $--$ & 0.5774 & 0.4804 & 0.6938 & Monotone plateau \\ \hline
		&0.25 & 0.66 & 5		  & 1.0560 & 1.0775   & 0.2702& 0.2027  & 0.3353 & Non-monotone basin \\
		\raisebox{1.6ex}[0pt]{Example 2} & 0.66 & 0 & 5			& 0.7545  & 2.5023 & 0.5774  & 0.4674 & 0.7130 & Non-monotone plateau \\ \hline
    &0.25 & 0.52 	& 5 		& 3.0723 & $--$		& 0.3246 	&  0.3109 & 0.3382 & Monotone basin \\
	\raisebox{1.6ex}[0pt]{Example 3} & 0.52 & 0 & 5 & 0.7545 & 0.4154 & 0.5774  & 0.4674 & 0.7130 & Non-monotone plateau \\ \hline
\end{tabular}
\end{table}


\subsection{Influence of the flux rate and the initial saturation}
In example 1, 2 and 3, the effect of the gravity has been neglected. In this section, we take the gravity into account to study the influences of flux rate and initial saturation. Since the smoothed monitor function has a better balance between the smooth and the steep regions than the arc-length monitor function, in the Example 4 and Example 5, we will only consider the smoothed monitor function with $N = 200$.

Consider $H(u)$ and $p_c(u)$ as in (\ref{eqn:symfh}) with the flux function $F(u)$ replaced by
\begin{align}
	F(u) = \frac{u^2}{u^2 + M (1 - u)^2}[ v_T + C(1 - u^2) ],
\end{align}
where $\epsilon = 10^{-3}$, $M = 10$, and $C = 10$ is a positive constant that accounts for gravity.
The initial condition is
\begin{align}
	\label{eqn:fluxinitc}
	u(x, 0) = u_0 + 0.5 (u_B - u_0) (1.0 - \tanh(200 x)).
\end{align}

\textbf{Example 4.} $ \tau = 3.3812, v_T = 1.0, 0.6, 0.4, 0.1, u_0 = 0, T = 1$.

In Ref. \cite{dicarlo2004experimental}, DiCarlo carried out a series of experiments with different infiltrating fluxes. At the highest ($2.0\times10^{-3} \mathrm{[m s^{-1}]}$) and lowest ($1.32\times10^{-7} \mathrm{[m s^{-1}]}$) fluxes, the profiles are monotonic and no saturation overshoot is observed. For the intermediate fluxes which exhibit saturation overshoot, as flux decreases both the tip and tail saturations decrease continuously. In this example, we use the simplified gravity model to show how solution varies with the change of flux.

\begin{figure}[!htbp]
\begin{center}
    {\includegraphics[width=2.4in] {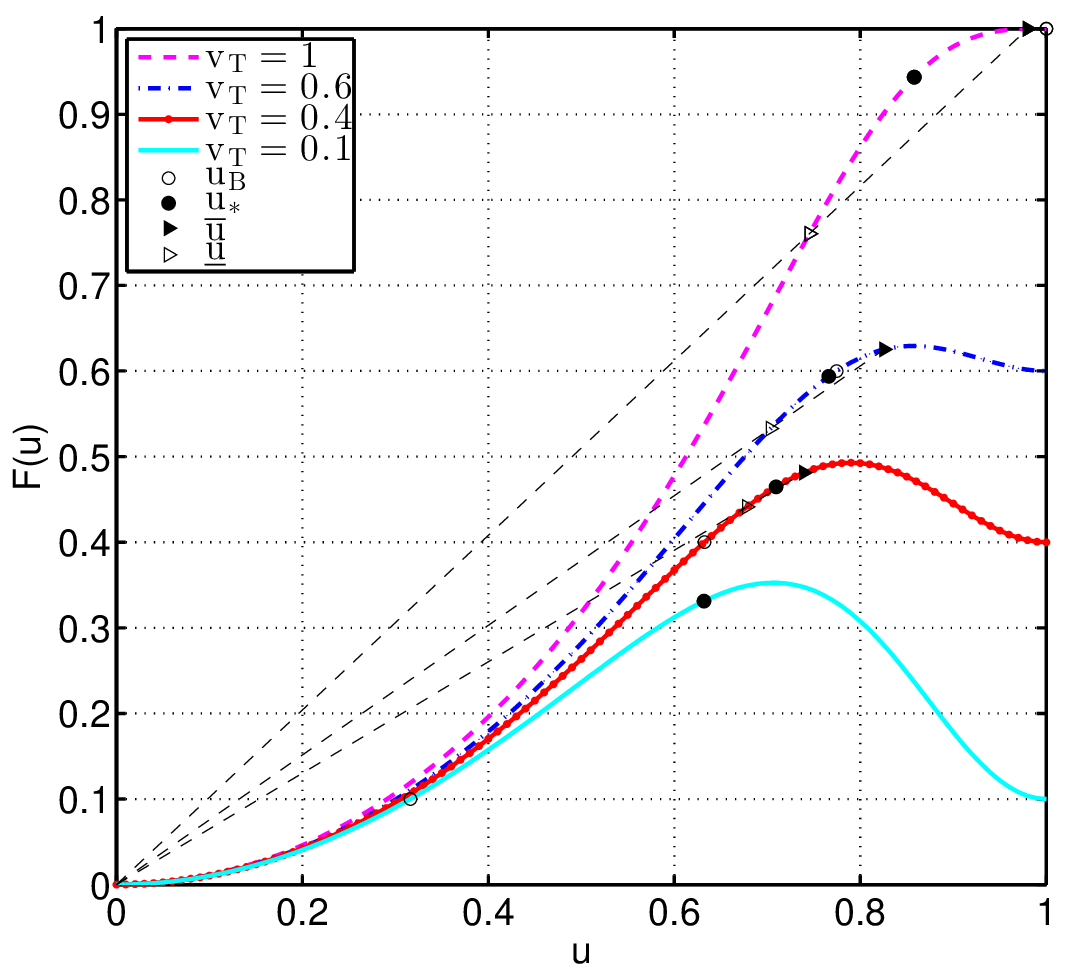}}
	{\includegraphics[width=2.4in] {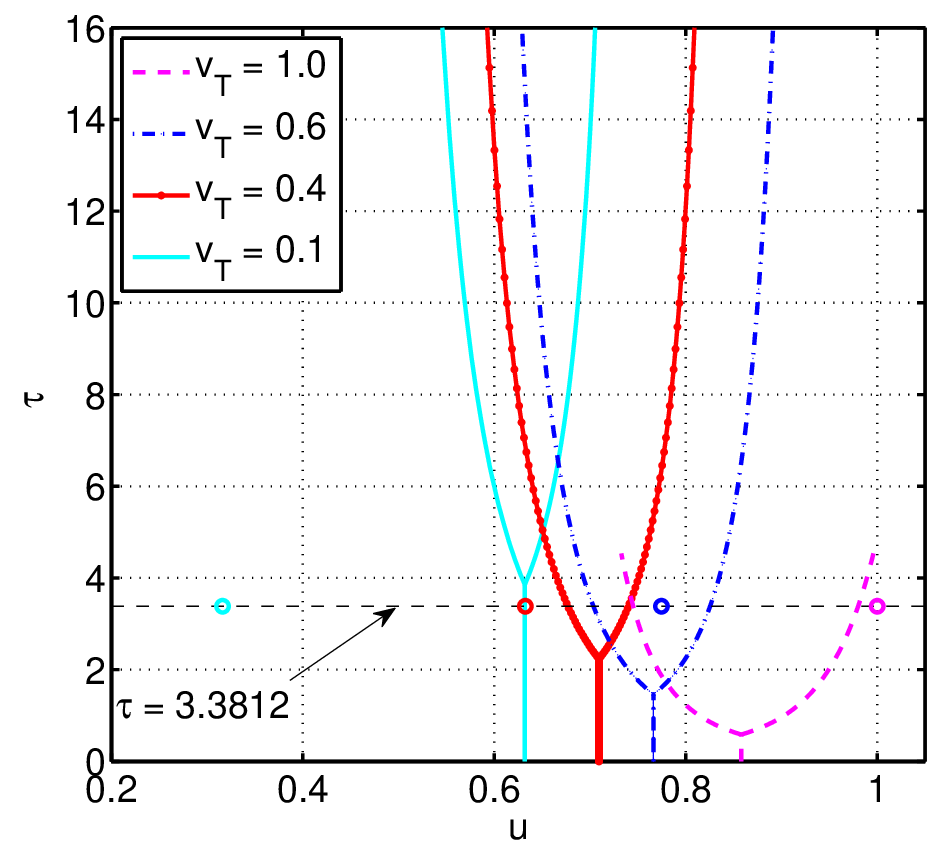}}
    {\includegraphics[width=2.4in] {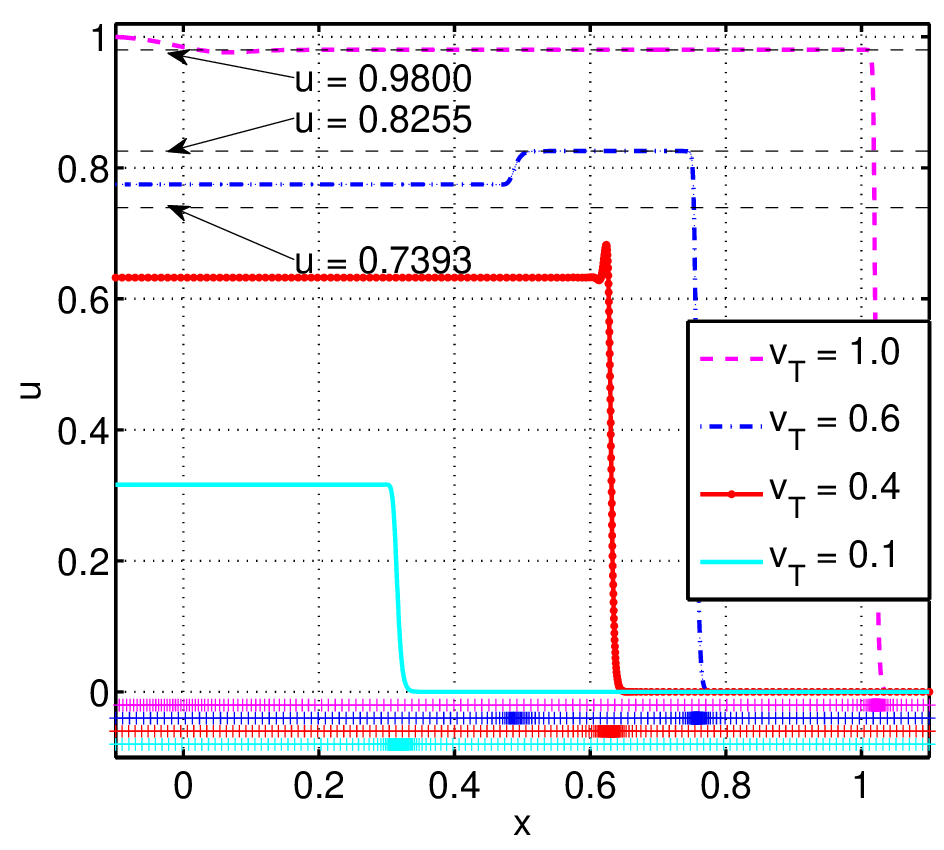}}
	 \caption{Flux functions (left), bifurcation diagrams (center) and numerical solutions (right) for different values of $v_T$.  \label{fig:fluxdiagcomp}}
\end{center}
\end{figure}

Next we fix the initial saturation $u_0 = 0$ and consider different values of $v_T$. The flux function $F(u)$ corresponding to $v_T$ are plotted in Fig. \ref{fig:fluxdiagcomp} (left). As can be seen, when $v_T = 1.0$, the flux function is strictly increasing, when $v_T$ is smaller, the flux function becomes non-monotonic, and the value of $u_\alpha$ becomes lower and results in different bifurcation diagrams. Fig. \ref{fig:fluxdiagcomp} (center) shows that the value of $\tau_*$ increases with decreasing $v_T$. By solving $F(u) = v_T$, the boundary saturation $u_B$ corresponding to different $v_T$ can be defined. For $\tau = 3.3812$, the TW results are given in Table \ref{tab:fluxtauu}.

When $v_T = 1$, we have $\tau > \tau_*$ and $u_B > \bar{u}(\tau)$. Thus $u_B$ and $u_0$ are connected by a rarefaction wave from $u_B$ down to $\bar{u}$ trailing an undercompressive shock from $\bar{u}$ down to $u_0$. When $v_T = 0.6$, we have $\underline{u}(\tau) < u_B < \bar{u}(\tau)$, and therefore a non-monotone plateau of height $\bar{u} = 0.8255$ exists. When $v_T = 0.4$, we have $u_B  < \underline{u}$ and $\tau_s < \tau < \tau_*$, hence $u_B$ is a spiral point: there exists an admissible Lax shock connecting $u_B$ and $u_0$ with oscillations near $u_B$. When $v_T = 0.1$, $u_B < u_\alpha$ and $\tau < \tau_*$, this combination results in a monotonic profile. The numerical solutions are plotted in Fig. \ref{fig:fluxdiagcomp} (right). The computed profiles have good agreement with the TW results.

\textbf{Example 5.} $\tau = 2.13, v_T = 0.6, u_B = 0.75, u_0 = 0, 0.1, 0.2, 0.25$, $T = 1$.

The laboratory experiments in Ref. \cite{dicarlo2004experimental} also show that saturation overshoot decreases quickly with increasing initial water saturation. In this example, we study the influence of the initial saturation $u_0$.

\begin{figure}[!htbp]
\begin{center}
   {\includegraphics[width=2.7in] {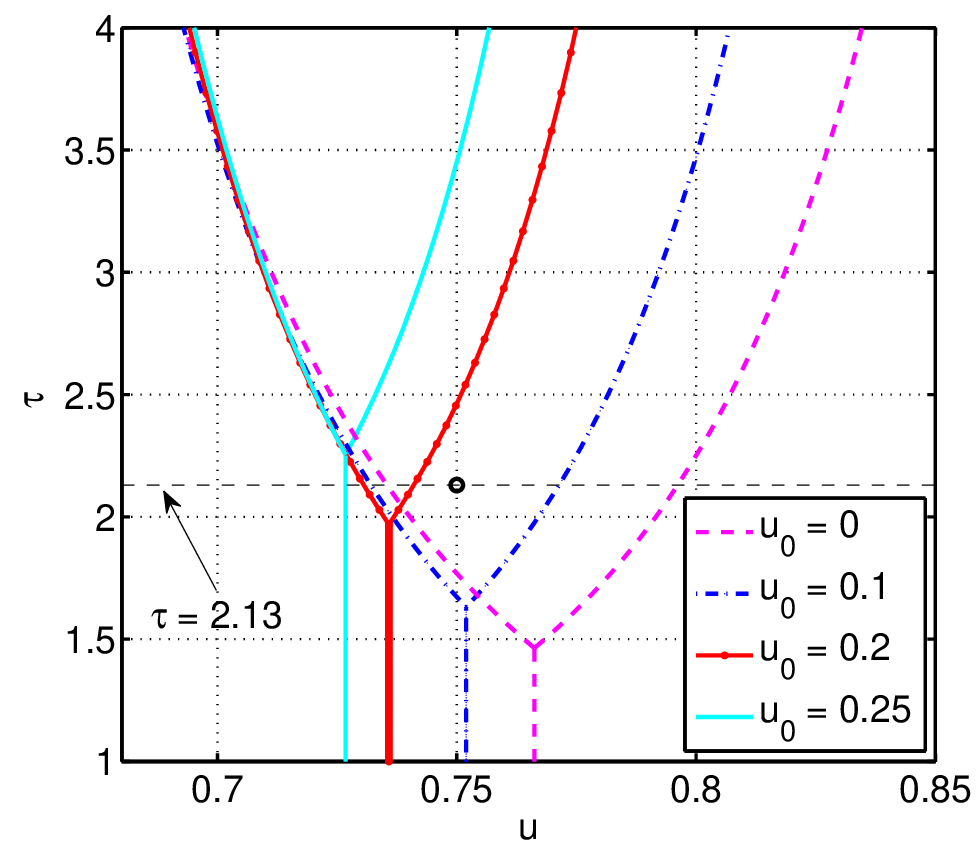}}\quad \quad
   {\includegraphics[width=2.7in] {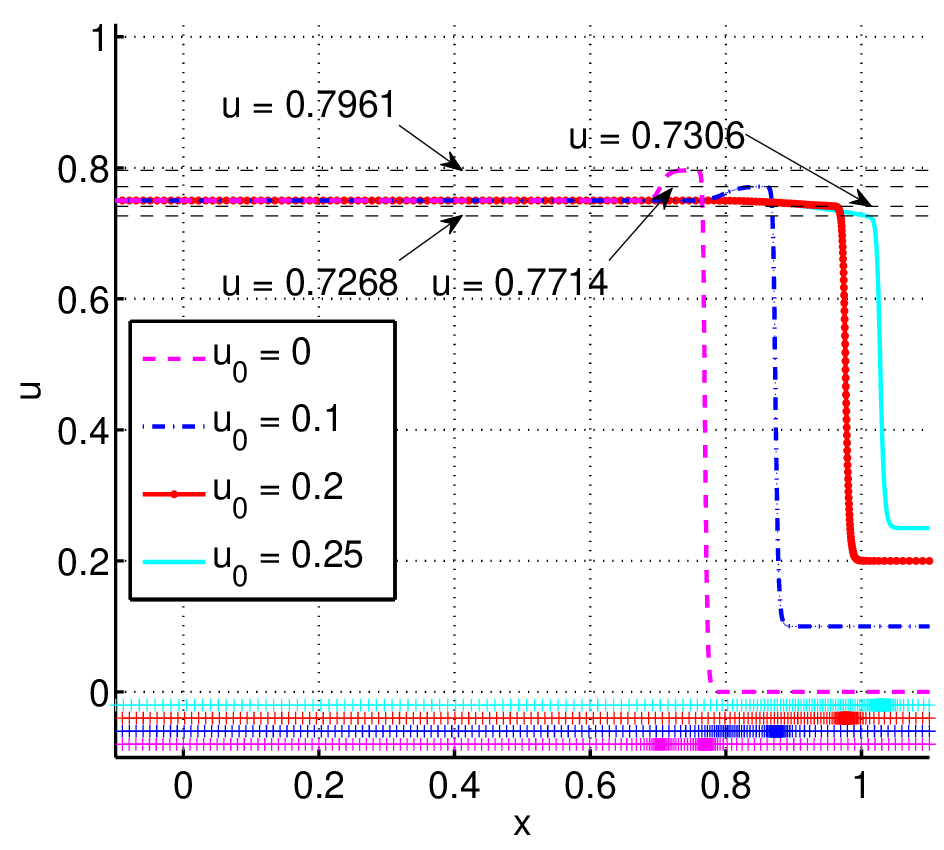}}
	 \caption{Bifurcation diagrams (left) and numerical solutions (right) for different values of $u_0$.\label{fig:initdiagcomp}}
\end{center}
\end{figure}

In Fig. \ref{fig:initdiagcomp} (left) the bifurcation diagrams for various values of $u_0$ are presented. Choosing $u_B = 0.75, \tau = 2.13$, the TW results are presented in Table \ref{tab:inittauu}. When $u_0 = 0$ and $0.1$, since $\underline{u} < u_B < \bar{u}$ and $\tau > \tau_*$, $u_B$ and $u_0$ are connected by non-monotone plateaus. The plateau height of $u_0 = 0$ is higher than that of $u_0 = 0.1$. When $u_0 = 0.2$, we have $u_B > \bar{u}$, therefore there exists a rarefaction wave from $u_B$ down to $\bar{u}$ trailing an undercompressive shock from $\bar{u}$ down to $u_0$. The plateau height $\bar{u}$ is lower than $u_0 = 0$ and $0.1$. When $u_r = 0.25$, since $\tau = 2.13$ is smaller than $\tau_* = 2.2537$, we have an admissible Lax shock connecting $u_B$ and $u_0$. The turning point is $u_\alpha = 0.7268$.

In Fig. \ref{fig:initdiagcomp} (right), we present the numerical solutions. As can be seen, the computed plateau heights and the turning points agree well with those in Table \ref{tab:inittauu}. The grids plotted in the lower part of the figure clearly indicate the critical regions of the solutions.

\begin{table}
	\caption{\label{tab:fluxtauu}Travelling wave results for different values of $v_T$ with $\tau = 3.3812, u_0 = 0$.}
  \center
  \begin{tabular}{|c|cccccc|c|}\hline
		$v_T$ & $u_B$ & $u_\alpha$ & $\tau_*$ & $\tau_s$ & $\underline{u}(\tau)$ & $\bar{u}(\tau)$ & Wave description\\ \hline
		$1.0$	& $1.0000$	& $0.8580$ 	& $0.5848$ & $--$
& $0.7452$ & $0.9800$ & Monotone plateau \\
		$0.6$ & $0.7746$  & $0.7662$ & $1.4633$  & $2.3406$
& $0.7035$ & $0.8255$ & Non-monotone plateau \\
		$0.4$ & $0.6325$  & $0.7093$ & $2.2330$ & $1.1512$
& $0.6779$ & $0.7393$ & Non-monotone overshoot \\
		$0.1$ & $0.3162$  & $0.6318$ & $3.8537$ & $2.6097$ & $--$ & $--$ & Monotone, no plateau
\\ \hline
\end{tabular}
\end{table}

\begin{table}
	\caption{\label{tab:inittauu}Travelling wave results for different values of $u_0$ with $\tau = 2.13, u_B = 0.75$.}
  \center
  \begin{tabular}{|c|ccccc|c|}\hline
		$u_0$  & $u_\alpha$ & $\tau_*$  & $\tau_s$ &
$\underline{u}(\tau)$ & $\bar{u}(\tau)$ & Wave description \\ \hline
		$0.00$ & $0.7662$	& $1.4633$ & $2.0371$
& $0.7354$ & $0.7961$ & Non-monotone plateau \\
		$0.10$ & $0.7519$ & $1.6312$ & $4.0993$
& $0.7320$ & $0.7714$ & Non-monotone plateau \\
		$0.20$ & $0.7358$ & $1.9688$ & $--$
& $0.7306$ & $0.7410$ & Monotone plateau \\
		$0.25$ & $0.7268$ & $2.2537$ & $--$
& $--$ & $--$ & Monotone, no plateau
\\ \hline
\end{tabular}
\end{table}

\subsection{Numerical solutions of the MBL equation using the Brooks-Corey model}

Ref. \cite{dicarlo2004experimental} presented snapshots of the saturation profiles for different fluxes infiltrating into initially dry $20/30$ sand. It's observed that for the highest ($v_T = 2.0\times{10^{-3}} \mathrm{~[m~s^{-1}]}$) and lowest ($v_T = 1.32\times{10^{-7}}\mathrm{~[m~s^{-1}]}$) fluxes the saturation profiles are monotonic with distance and no saturation overshoot is observed, while all of the intermediate fluxes ($v_T = 1.32\times 10^{-3}\mathrm{~[m~s^{-1}]}$, $1.32\times 10^{-4}\mathrm{~[m~s^{-1}]}$, $1.32\times 10^{-5}\mathrm{~[m~s^{-1}]}$, $1.32\times 10^{-6}\mathrm{~[m~s^{-1}]}$) exhibit saturation overshoots. In this section we study the numerical solutions of (\ref{eqn:fullmbl}) with flux function (\ref{eqn:flux}) and diffusion function (\ref{eqn:diff}). The physical parameters of the $20/30$ sand \cite{dicarlo2004experimental,schroth1996characterization} as well as the constants and the Brooks-Corey model \cite{brooks1966properties} are listed in Table \ref{tab:2030sand} and Table \ref{tab:model}.

\textbf{Example 6.} Full equation using the Brooks-Corey type model.

\begin{figure}[!htbp]
\begin{center}
		{\includegraphics[width=2.7in] {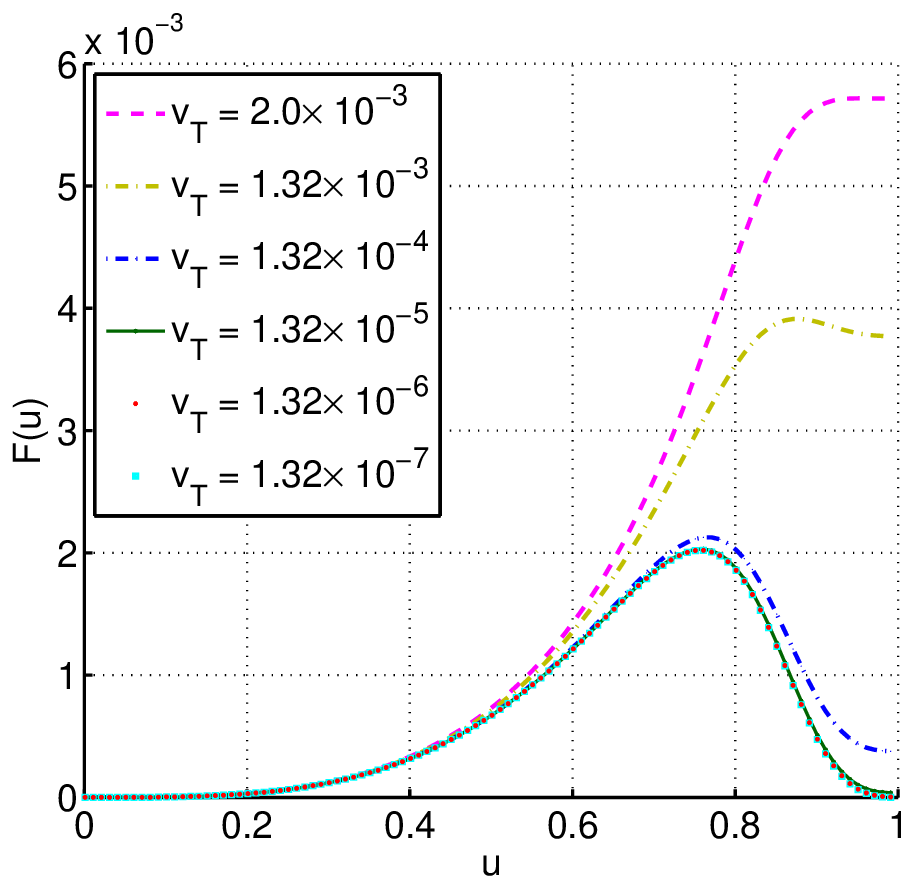}}\quad \quad
   {\includegraphics[width=2.7in] {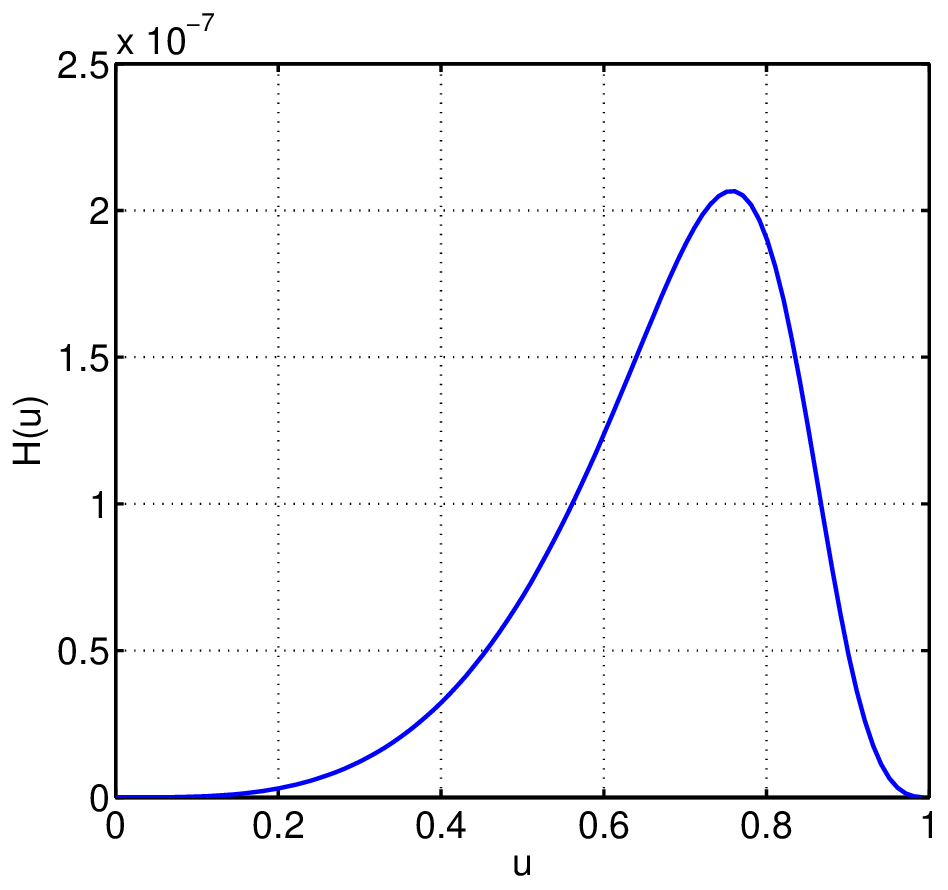}}
	 \caption{The graphs of $F(u)$ for different values of $v_T$ (left) and $H(u)$ (right).
 \label{fig:fullmodelfh}}
\end{center}
\end{figure}

Using the physical parameters and the Brooks-Corey type model in Table \ref{tab:2030sand} and Table \ref{tab:model}, the flux $F(u)$ with different values of $v_T$ and the diffusion function $H(u)$ are plotted in Fig. \ref{fig:fullmodelfh}. The degeneracy of $H(u)$ at $u = 0$ makes (\ref{eqn:fullmbl}) difficult to solve. From Fig. 1 in Ref. \cite{dicarlo2007capillary} we get the initial capillary pressure $p_c^{0} \approx 1600 \text{[Pa]}$, when initial water saturation $u_0 = 0.003$, using parameters in the imbibition process we get the Brooks-Corey capillary pressure $p_{c}(u_0) = 1566 \text{[Pa]}$. Thus in the numerical simulations the initial saturation $u_0 = 0.003$ is adopted. This initial saturation is also the measured value in Ref. \cite{yao2001stability}.

\begin{table}
  \caption{\label{tab:2030sand}Physical parameters for $20/30$ sand.}
  \center
  \begin{tabular}{|cccccc|ccc|}\hline
  {}& {} & {} & \multicolumn{3}{c|}{Drainage}&\multicolumn{3}{c|}{Imbibition}\\\cline{4-9}
{Sand} & $\kappa~\mathrm{[m~s^{-1}]}$ & $\phi$~[-] & $u_{re}$~[-] & $\lambda$~[-]& $p_d ~\text{[Pa]}$ & $u_{re}$~[-] & $\lambda$~[-] & $p_d ~\text{[Pa]}$\\\hline
    $20/30$ & $2.5\times 10^{-3}$ & $0.35$ & $0$ & $5.57$ & $850$ & $0$ & $5$ & $490$ \\\hline
  \end{tabular}
\end{table}

\begin{table}
  \caption{\label{tab:model}Constants and the Brooks-Corey model.}
  \center
  \begin{tabular}{|c|cc|}\hline
  Density $\mathrm{[kg~m^{-3}]}$ & $\rho_w = 998.21$ & $\rho_n = 1.2754$  \\
  Viscosity $\mathrm{[kg~m^{-1} s^{-1}]}$ & $\mu_w = 1.002\times 10^{-3}$ & $\mu_n = 1.82\times 10^{-5}$  \\
  Mobility $\mathrm{[m~s~kg^{-1}]}$ & $\lambda_w = \frac{K k_{rw}}{\mu_w}$ & $\lambda_n = \frac{K k_{rn}}{\mu_n}$ \\
  Constants & $g = 9.81 ~\mathrm{[m~s^{-2}]}$  & $K = \frac{\kappa \mu_w}{\rho_w g}$ $~\mathrm{[m^2]}$\\\hline
   & Capillary pressure  &{Relative permeability } \\\hline
   & $u_e = \frac{u - u_{re}}{1 - u_{re}}$ &  $k_{rw} = u_e^{\frac{2+ 3 \lambda}{\lambda}}$ \ \\
  \raisebox{1.6ex}[0pt]{Brooks-Corey model} & {
  $p_c = p_d u_e^{-\frac{1}{\lambda}},~~\mathrm{for}~p_c > p_d$}& $k_{rn} = (1-u_e)^2 (1-u_e^{\frac{2+\lambda}{\lambda}})$ \\\hline
  \end{tabular}
\end{table}

\begin{figure}[!htbp]
\begin{center}
		{\includegraphics[width=2.7in] {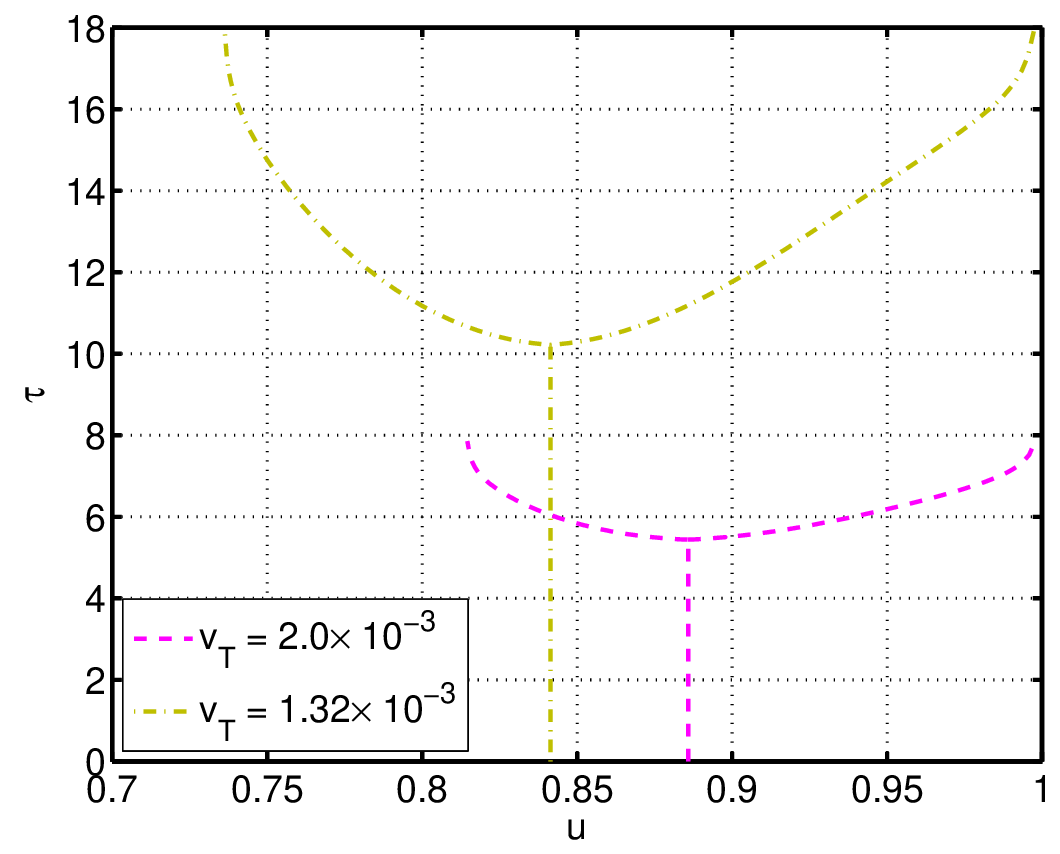}}\quad \quad
   {\includegraphics[width=2.7in] {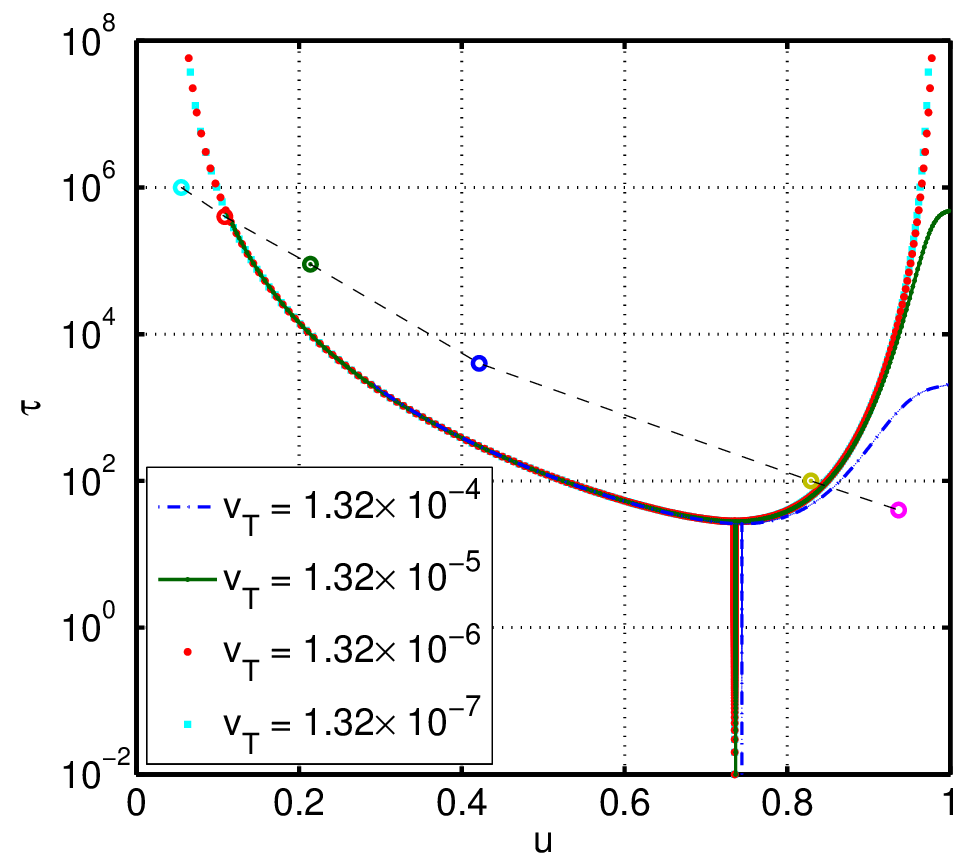}}
	 \caption{The bifurcation diagrams for different values of $v_T$.
 \label{fig:fullmodeldiag123456}}
\end{center}
\end{figure}

\begin{figure}[!htbp]
\begin{center}
   {\includegraphics[width=2.7in] {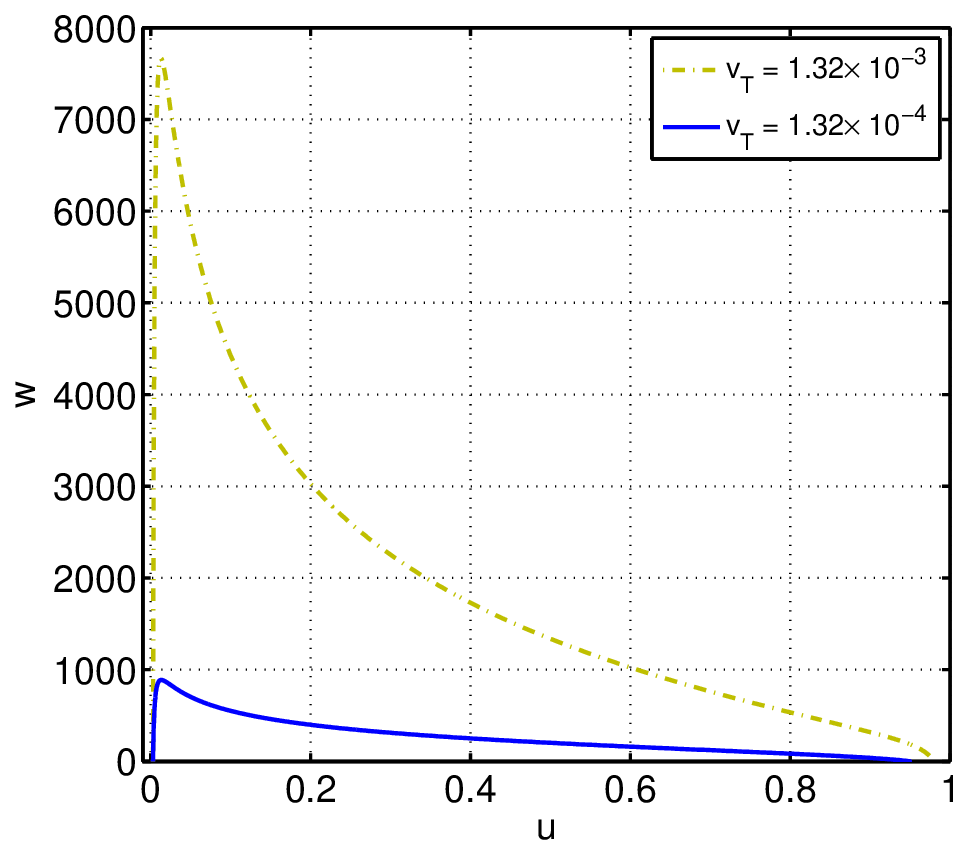}}
   {\includegraphics[width=2.7in] {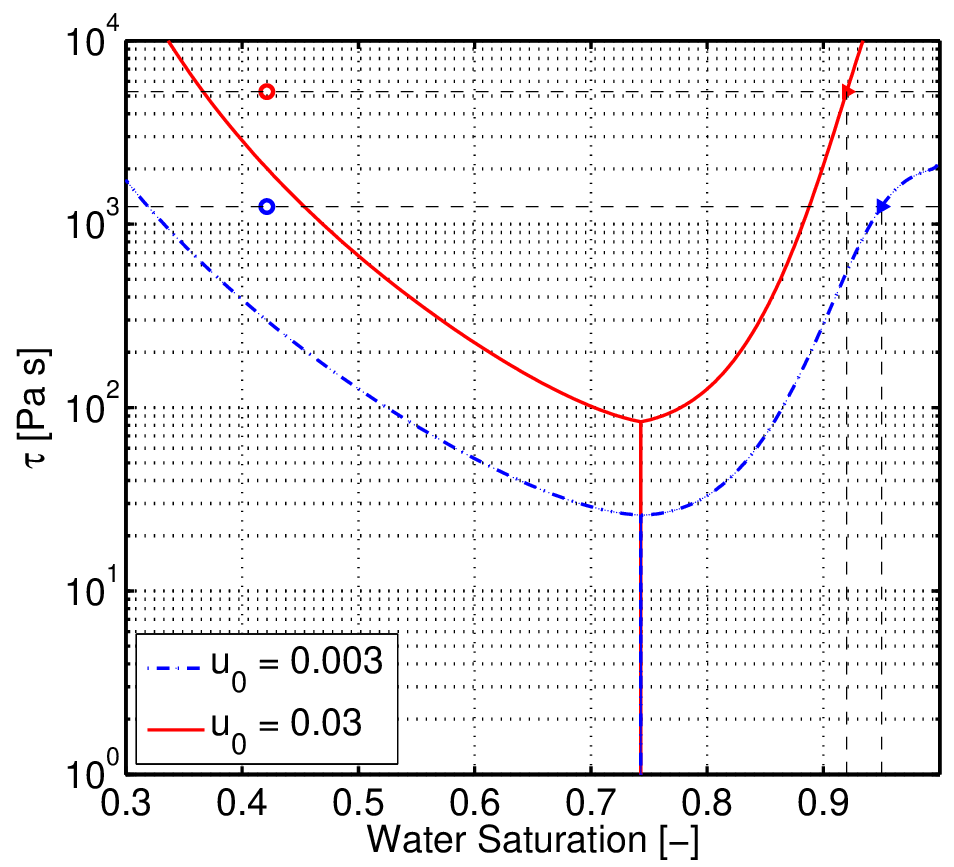}}
	 \caption{Left: graphs of $w(u)$ with $v_T = 1.32\times 10^{-3}$ and $1.32\times 10^{-4}$; Right: the bifurcation diagrams for initial saturation $u_0 = 0.003, 0.03$ with $v_T = 1.32\times 10^{-4}$.
 \label{fig:fullmodeldiagurinitwu}}
\end{center}
\end{figure}

Fig. \ref{fig:fullmodelfh} (left) shows that, the flux function $F(u)$ differs a lot for $v_T = 2.0\times 10^{-3}, 1.32\times10^{-3}, 1.32\times10^{-4}$. Therefore, the bifurcation diagrams for these three cases are different from each other, see in Fig. \ref{fig:fullmodeldiag123456}. In order to simulate the saturation overshoot phenomenon, the dynamic coefficient $\tau$ has to be determined. We plot the $\tau - u$ pairs used in Ref. \cite{zhang2016mimetic}. Notice that the $\tau\text{-} u$ pairs are beyond the scope of the bifurcation diagrams for $v_T = 2.0\times10^{-3}, 1.32\times 10^{-3}$ and $1.32\times10^{-4}$. Thus we have to choose new values of $\tau$ according to the bifurcation diagrams. For $v_T = 1.32\times 10^{-3}, 1.32\times 10^{-4}$, let the corresponding plateau saturations $\bar{u}$ equal $0.98$ and $0.95$, by solving (\ref{eqn:slope}) and (\ref{eqn:slopeboundary}), the values of $\tau$ are $15.81$ and $1246$, respectively. The graphs of $w(u)$ are presented in Fig. \ref{fig:fullmodeldiagurinitwu} (left). Notice $w = -u_x$, such high values of $-u_x$ for $v_T = 1.32\times 10^{-3}$ can't be achieved using moving mesh method with $N \leq 800$ or uniform grid with $N \leq 4000$. In the following we will only consider $v_T = 1.32\times{10^{-4}}$.
The bifurcations diagrams for $v_T = 1.32\times 10^{-4}$ with two different initial saturations $u_0 = 0.003$ and $0.03$ are presented in Fig. \ref{fig:fullmodeldiagurinitwu} (right).  
Assuming that after a long time, the saturation at the left boundary reaches the equilibrium state, we set $\frac{\partial u}{\partial t}|_{x_l} = 0, \frac{\partial u}{\partial x}|_{x_l} = 0$. Then we obtain
\begin{align}\label{eqn:boundsatu}
  v_T f(u) + \lambda_n(u) f(u) (\rho_w - \rho_n)g = v_T.
\end{align}
Solving (\ref{eqn:boundsatu}) we get the boundary saturation $u_B$ corresponding to $v_T$.

\begin{figure}[!htbp]
\begin{center}
    {\includegraphics[width=2.7in] {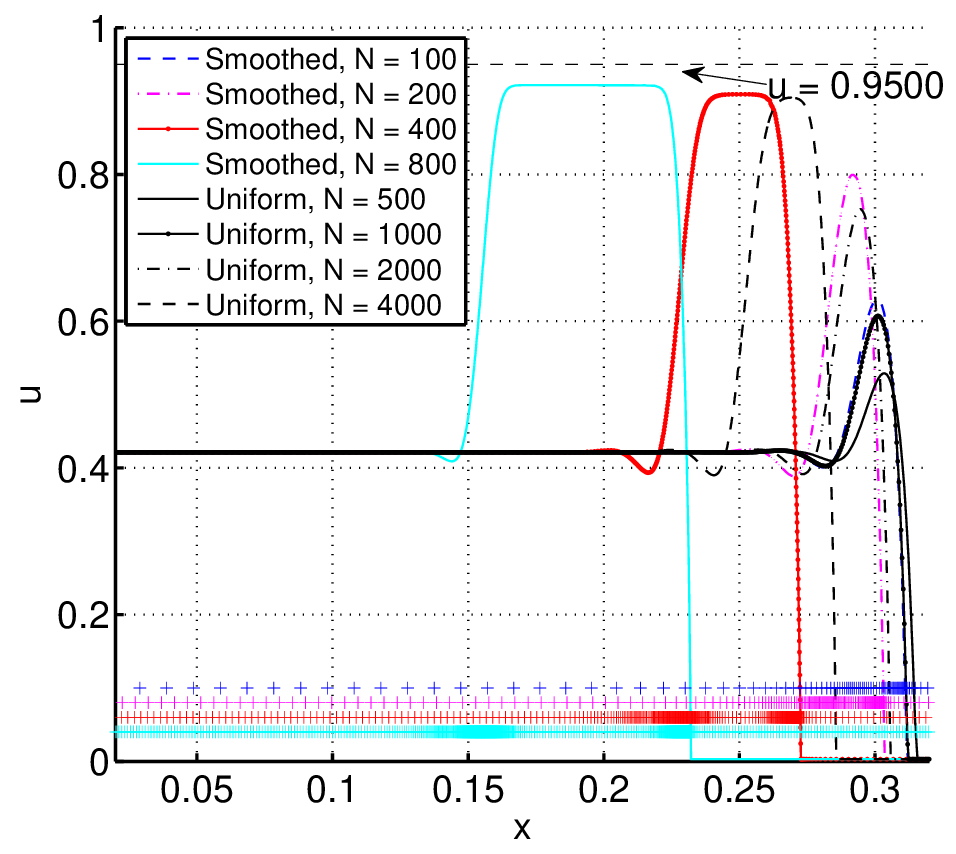}} \quad \quad
	 {\includegraphics[width=2.7in] {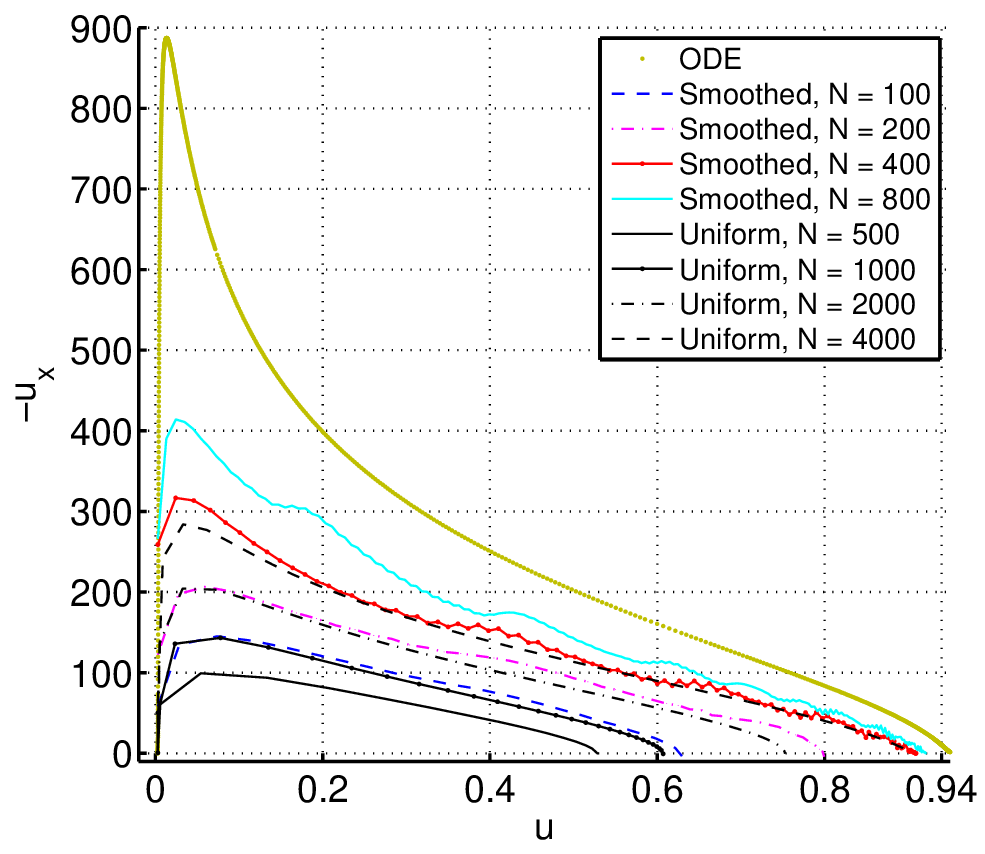}}
   {\includegraphics[width=2.7in] {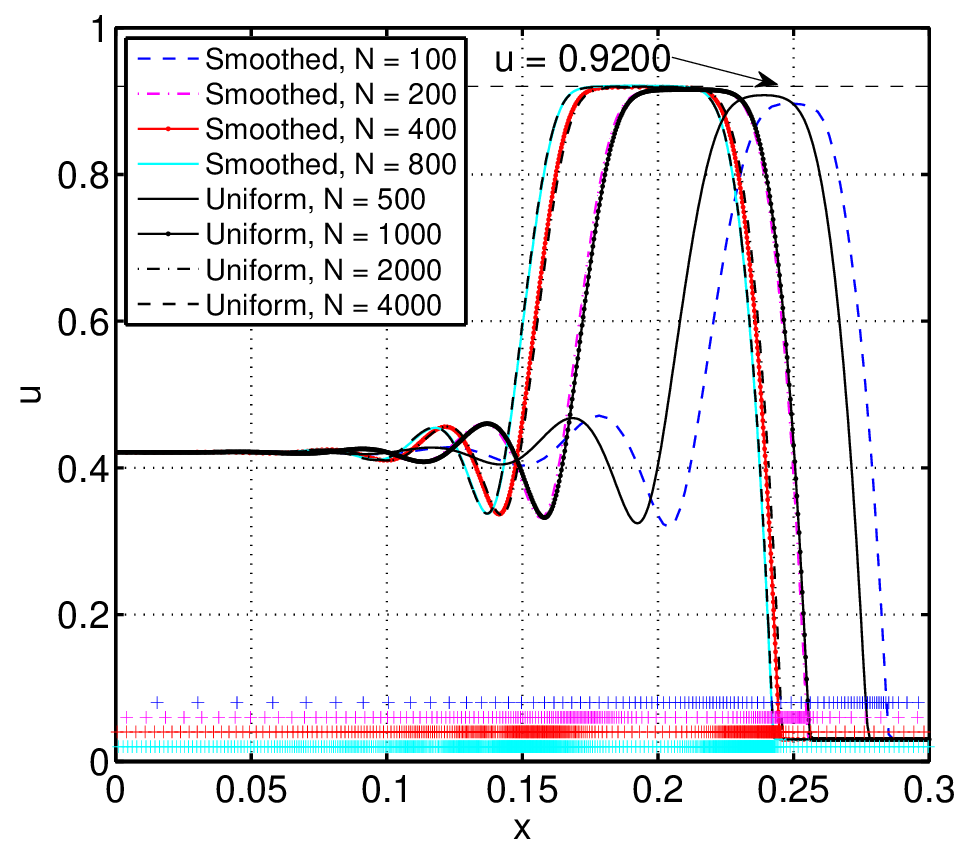}}\quad \quad
   {\includegraphics[width=2.7in] {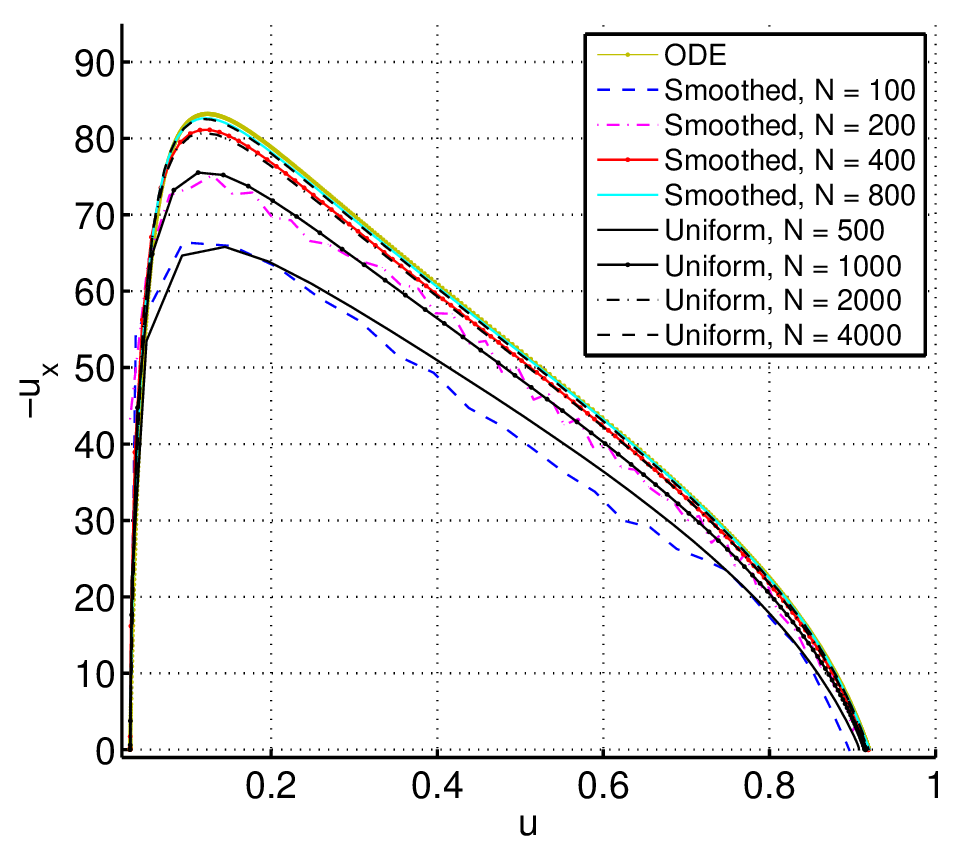}}
	\caption{Example 6: solutions (left) and values of $-u_x$ at the right boundary of the undercompressive shock (right) of the full model with $v_T = 1.32\times10^{-4}$. $u_0 = 0.003$ (top), $u_0 = 0.03$ (bottom).
 \label{fig:case4ur003ur03}}
\end{center}
\end{figure}

Using the initial condition
\begin{align}
	\label{eqn:fullmodelinitc}
	u(x, 0) = u_0 + 0.5 (u_B - u_0) (1.0 - \tanh(200 x)), \quad x\in[-0.05, 0.35],
\end{align}
the TW results obtained from Section \ref{sec:tw} are listed in Table \ref{tab:twcase4ur003ur03}, the numerical results at $t = 350 [s]$ are shown in Fig. \ref{fig:case4ur003ur03}. When $u_0 = 0.003$, it shows the more grid points we use, the more accurate is the plateau height. The uniform grid with $N = 4000$ only gives an overshoot which is lower than that of the moving mesh method with $N = 400$. The plateau height obtained by the moving mesh method with $N = 800$ is also lower than $\bar{u}$ in Table \ref{tab:twcase4ur003ur03}. This may be caused by the lack of grid points near the undercompressive shock: there are more grid points near the undercompressive shock for the moving mesh method than for the uniform grid, Fig. \ref{fig:case4ur003ur03} (right) shows that the values of $-u_x$ obtained by the moving mesh method are higher than those using the uniform grids.

\begin{figure}[!htbp]
\begin{center}
    {\includegraphics[width=2.7in] {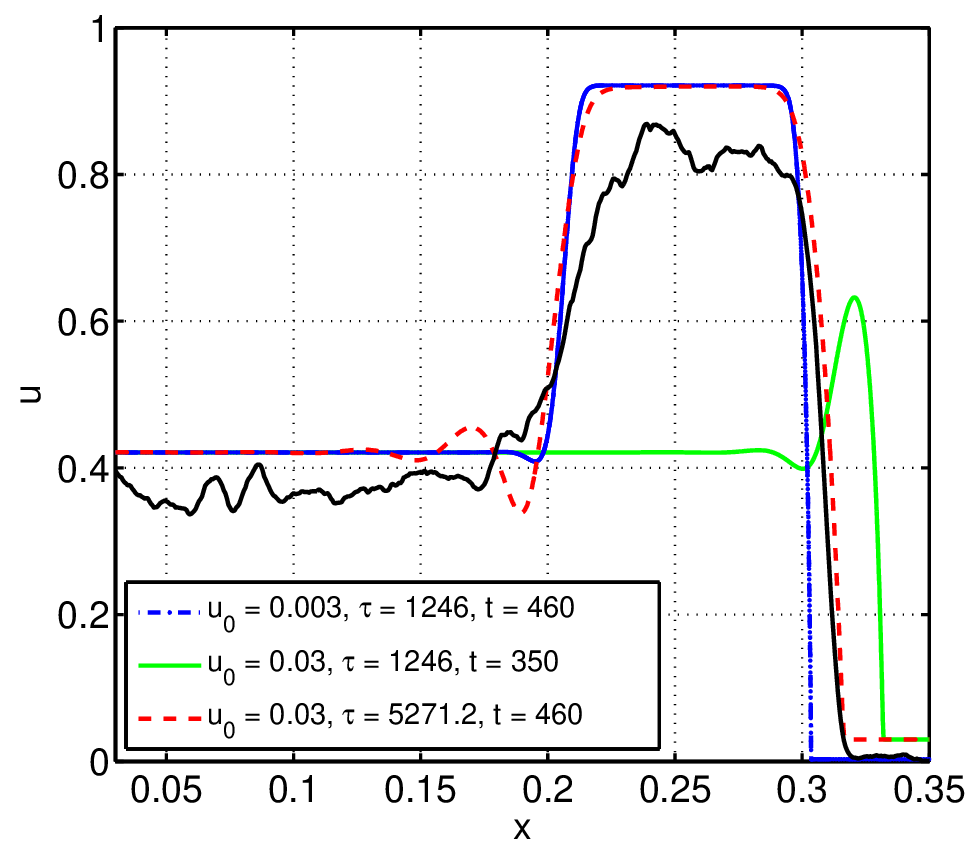}}
	\caption{Comparisons between experimental result and numerical solutions obtained for $(u_0, \tau, t) = (0.003, 1246, 460)$, $(0.03, 1246, 350)$, $(0.03, 5271, 460)$ using moving mesh method with $N = 800$.
 \label{fig:case4plots3}}
\end{center}
\end{figure}

In Fig. \ref{fig:case4plots3} we compare the saturation profiles obtained by different values of $(u_0, \tau, t) = (0.003, 1246, 460)$, $(0.03, 1246, 350)$, $(0.03, 5271, 460)$ with the experimental results from Ref. \cite{dicarlo2004experimental}. When $u_0 = 0.03$ and $\tau = 1246$, Table \ref{tab:twcase4ur003ur03} shows $u_B < \underline{u}$, $\tau > \tau_s$, thus instead of a plateau, the dynamic coefficient $\tau = 1246$ only gives an overshoot. When $u_0 = 0.03$ is fixed, a plateau is obtained using a higher value of $\tau = 5271$. For $(u_0, \tau, t) = (0.03, 1246, 350)$ and $(0.03, 5271, 460)$, since $\tau > \tau_s$, oscillations appear near $u_B = 0.4212$. When $u_0 = 0.003$, the computed plateau height is $\bar{u}_c = 0.9214$, thus $\tau > \tau_{sc} = 614.5$, a small oscillation still exists near $u_B$. The plateau and tail saturations obtained with $(u_0, \tau, t) = (0.003, 1246, 460)$ are higher than the experimental values, we attribute this to the limitation of the Brooks-Corey model. Note that, the end time for $(u_0, \tau, t) = (0.003, 1246, 460)$ is very near to the calculated end time $t = 409.7$ in Ref. \cite{zhang2016mimetic}, this validates the accuracy of our approach.

\begin{table}
  \caption{\label{tab:twcase4ur003ur03}Travelling wave results for $u_0 = 0.003, 0.03$ with $v_T = 1.32\times 10^{-4}$}
  \center
  \begin{tabular}{|c|ccccccc|c|}\hline
	$u_0$  & $u_B$ & $\tau$ & $\tau_*$ & $\tau_s$ & $u_*$  & $\underline{u}$ & $\bar{u}$ & Wave description \\ \hline
    0.003 & 0.4212 & 1246 & 25.89 & 2117 & 0.7440 & 0.3190 & 0.9500 & Non-monotone plateau \\
    0.03  & 0.4212 & 1246  & 83.52 & 368.4 & 0.7429 & 0.4533 & 0.8878 & Non-monotone overshoot \\
    0.03 & 0.4212  & 5271 & 83.52 & 588.9 & 0.7429 & 0.3664 & 0.9200 & Non-monotone plateau
 \\ \hline
\end{tabular}
\end{table}

\section{Conclusions}
In this paper we applied an adaptive moving mesh technique to solve the two-phase flow equation incorporating the dynamic capillary pressure term. The moving mesh method successfully captured the monotone and non-monotone solutions with high accuracy. Comparisons between numerical results show that to achieve the same accuracy, the moving mesh method needs approximately a factor of 5-10 fewer grid points than the uniform case. The computed saturation profiles and grid trajectories also illustrate different features of the smoothed monitor function and the arc-length monitor function. The arc-length monitor function have higher accuracy in steep regions, while the smoothed monitor function gives a better balance between the smooth and the steep regions.

\section*{Acknowledgements}
The work of H. Zhang was supported by the China Scholarship Council (No. 201503170430). We would like to thank Prof. Iuliu Sorin Pop for providing us his code for the computation of $\tau$ in Example 6.

\bibliographystyle{elsarticle-num}
\bibliography{reference}
%
\end{document}